\documentclass[12pt]{amsart}
\usepackage{latexsym}
\usepackage{amssymb}
\usepackage{amsxtra}
\usepackage{amsmath}
\usepackage{amsfonts}
\usepackage{amscd}
\usepackage{verbatim}
\usepackage{color}

\newcommand{\beq}{\begin{equation}}
\newcommand{\eeq}{\end{equation}}
\newcommand{\beqa}{\begin{eqnarray}}
\newcommand{\eeqa}{\end{eqnarray}}
\newcommand{\beaa}{\begin{eqnarray*}}
\newcommand{\ben}{\begin{eqnarray*}}
\newcommand{\eaa}{\end{eqnarray*}}
\newcommand{\een}{\end{eqnarray*}}

\newcommand \nc {\newcommand}

\newtheorem{theorem}{Theorem}[section]
\newtheorem{lemma}[theorem]{Lemma}
\newtheorem{proposition}[theorem]{Proposition}
\newtheorem{corollary}[theorem]{Corollary}

\newtheorem{remark}[theorem]{Remark}
\newtheorem{conjecture}[theorem]{Conjecture}

\nc \thref{Theorem \ref}
\nc \leref{Lemma \ref}
\nc \prref{Proposition \ref}
\nc \coref{Corollary \ref}
\nc \deref{Definition \ref}
\nc \exref{Example \ref}
\nc \reref{Remark \ref}

\newcommand{\leftexp}[2]{{\vphantom{#2}}^{{\rm #1}}{#2}}
\newcommand{\leftbase}[2]{{\vphantom{#2}}_{{\rm #1}}{#2}}

\newcommand{\A}{\mathcal{A}}
\newcommand{\B}{\mathcal{B}}
\newcommand{\C}{\mathbb{C}}
\newcommand{\D}{\mathcal{D}}

\newcommand{\F}{\mathcal{F}}
\renewcommand{\H}{\mathcal{H}}

\newcommand{\M}{\mathcal{M}}

\renewcommand{\O}{\mathcal{O}}
\renewcommand{\P}{\mathbb{P}}
\newcommand{\QQ}{\mathbb{Q}}

\renewcommand{\S}{\mathcal{S}}
\newcommand{\T}{\mathcal{T}}

\newcommand{\Z}{\mathbb{Z}}

\newcommand{\f}{\mathbf{f}}

\newcommand{\q}{\mathbf{q}}

\def\d{\partial}

\def\iso{\cong}

\def\({\left(}
\def\){\right)}
\def\[{\left[}
\def\]{\right]}
\def\<{\left\langle}
\def\>{\right\rangle}

\def\one{{\bf 1}}
\def\gl{\lambda}
\def\la{\lambda}
\def\si{\sigma}
\def\Si{\Sigma}
\def\ge{\epsilon}
\def\ga{\alpha}

\begin{document}
\title[Orbifold $\P^1$ and quasi-modular forms]{Gromov-Witten theory of elliptic orbifold $\mathbb{P}^1$ and quasi-modular forms}
\author{Todor Milanov \& Yongbin Ruan}
\address{IPMU\\  Todai Institute for Advanced Studies (TODIAS)\\ University of Tokyo\\Kashiwa, Chiba 277-8583, Japan}
\email{todor.milanov@ipmu.jp}
\address{Department of Mathematics\\ University of Michigan\\ Ann Arbor,
MI}\email{ruan@umich.edu}

\maketitle
\tableofcontents
\section{Introduction}

A major problem in geometry and physics is to compute the Gromov-Witten invariants of a given target manifold.  In general, this is a complicated problem. However, in certain special situations, the computations lead to beautiful objects, such as modular forms. It is clearly an important problem to locate all these special examples where the modularity exists. The simplest example of this phenomenon is the genus-1, degree-$d$ invariants $n_{1,d}$ of an
elliptic curve $E$.  It is well-known that their generating function can be expressed in terms of the Dedekind $\eta$-function
$$
\exp\Big( -\sum_{d\geq 1}n_{1,d} q^d\Big)=q^{-1/24}\ \eta(q).
$$
One can say much more.  Let us introduce some notation. Let $X$ be a projective manifold and $\overline{\M}_{g,k}(X, \beta)$ be the moduli space of genus $g$ stable maps with $k$ markings and fundamental class $\beta$. Let $e_i$ be the evaluation map at the $i$-th marked point $x_i$ and $\psi_i$ be the first Chern class of the cotangent line bundle at $x_i$. Choose a basis $\phi_i$ of $H^*(X, \QQ)$ with $\phi_0=1$. The numerical GW invariants are defined by
$$
\langle\tau_{l_1}(\phi_{i_1}), \dots, \tau_{l_k}(\phi_{i_k})\rangle^X_{g,\beta}=
\int_{[\overline{\M}_{g,k}(X, \beta)]^{vir}} \prod_i (e^*_i\phi_i)\psi^{l_i}_1.
$$
The above invariant is zero unless
$$
\sum_i ({\rm deg}(\phi_i)+2l_i)=2(c_1(TX)(\beta)+(3-n)(g-1)+k).
$$
The advantage of Calabi-Yau manifolds, such as the elliptic curve $E$, is that $c_1(TX)=0$ and hence the dimension constraint is independent of $\beta$. For the elliptic curve $E$, the degree $\beta$ can be identified with a non-negative integer $d$. Then, it is natural to define
\beq\label{cor:E}
\langle\tau_{l_1}(\phi_{i_1}), \dots, \tau_{l_k}(\phi_{i_k})\rangle^E_{g}(q)=
\sum_{d\geq 0}
\langle\tau_{l_1}(\phi_{i_1}), \dots, \tau_{l_k}(\phi_{i_k})\rangle^E_{g,d}q^d.
\eeq
The genus-1 invariant $n_{1,d}$ from above corresponds to $\langle\ \rangle^E_1(q)$. By the dilaton and the divisor equations, the invariants with insertion
         $\tau_1(1), \tau_0(\phi_{-1})$ can be deduced from other invariants. Without loss of generality, we assume that
$\tau_l(\phi_i)\neq \tau_1(1), \tau_0(\phi_{-1})$.
Then, Okounkov-Pandharipande \cite{OP1, OP2, OP3} showed that the invariant \eqref{cor:E} converges to a quasi-modular form of $SL_2(\Z)$ with the change of variable $q=e^{2\pi i\tau}$. Together with a result of Krawitz-Shen \cite{KS}, we shall prove the modularity for  another class of examples, the elliptic
orbifold $\P^1$ with weights $(3,3,3), (2,4,4), (2,3,6)$. These orbifolds are the quotients of an elliptic curve $E$. Our methods however, are completely different from the methods of Okounkov--Pandharipande.

To state the theorem, choose a basis $\phi_i, i=-1, 0, \dots$ of $H^*_{CR}$ such that $\phi_{-1}$ is the divisor class and $\phi_0=1$. In the above cases, $c_1(TX)=0$ and we can define
\beq\label{cor:X}
\langle\tau_{l_1}(\phi_{i_1}), \dots, \tau_{l_k}(\phi_{i_k})\rangle^X_{g}(q)
\eeq
similarly.
The main result of the current paper is the following modularity theorem.
\begin{theorem}\label{t1}
Suppose that $\tau_l(\phi)\neq \tau_1(1),\tau_0(\phi_{-1})$ and $X$ is one of the three elliptic orbifolds $\P^1$ from above. For any multi-indices $l_j, i_j$, the GW invariant \eqref{cor:X} converges to
a quasi-modular form of an appropriate weight for a finite index subgroup $\Gamma$ of $SL_2(\Z)$ under the change of variables $q=e^{2\pi i \tau/3}$, $e^{2\pi i \tau/4},$
$e^{2\pi i \tau/6},$ respectively (see section \ref{PXJ} for the subgroups $\Gamma$ and the weights of the quasi-modular forms).
\end{theorem}

We would like to remark that if we include insertions of the form $\tau_1(1)$ then a similar statement holds.
In this case however, we need to perform a dilaton shift which amounts to taking linear combinations of
the above invariants.

The modular invariance has been at the center of recent physical developments of Gromov-Witten theory by Klemm and his collaborators \cite{ABK, GKMW}. Some of the key ideas such as anti-holomorphic completion were directly inspired by their work, for which the authors express their special thanks. There is a  work of similar flavor by Coates-Iritani on modularity of GW invariants of local $\P^2$ \cite{CI}. We are informed that Paul Johnson has an independent approach to the results in this paper. We thank them for interesting discussions. When this paper is finished, we notice a related paper
of Costello-Li where they constructed a B-model high genus theory of elliptic curve and obtained corresponding mirror symmetry \cite{CL}. Finally, Satake--Takahashi \cite{ST} established an isomorphism between the quantum cohomology of the above orbifold projective lines and the Milnor rings of the simple elliptic singularities, which is an important step in our main construction (although we do not make use of their results).

\subsection{Relation to the work of Krawitz-Shen}
There is a companion article by Krawitz-Shen \cite{KS}. Together, we completely solved all the problems regarding the GW theory and related
topics for the above three classes of orbifolds. The idea is from  the Landau-Ginzburg/Calabi-Yau correspondence.
Since the general philosophy applies to many other examples, let us briefly outline it.

Recall that a polynomial $W$ is called {\em quasi-homogeneous}  if there are rational numbers $q_i$,
called the {\em degrees} or the {\em charges} of $x_i$, such that
$$
W(\lambda^{q_0}x_0,\lambda^{q_1}x_1, \dots, \lambda^{q_N}x_N)=\lambda W(x_0,x_1,\dots, x_N)
$$
for all $\lambda\in \C^*$. The polynomial $W$ is called {\em non-degenerate} if: (1) $W$
defines a unique singularity at zero; (2) the choice of $q_i$
is unique.  A diagonal matrix ${\rm diag}(\la_0,\lambda_1, \dots, \lambda_N)$ is called an {\em abelian} or
{\em diagonal symmetry} of $W$ if
\ben
W(\la_0x_0,\la_1 x_1, \dots, \la_N x_N)=W(x_0, x_1, \dots, x_N).
\een
The diagonal symmetries form a group $G_{max}$ which is always nontrivial since it contains the element
$$
J_W={\rm diag}(e^{2\pi i q_0}, e^{2\pi i q_1}, \dots, e^{2\pi i q_N}).
$$
When $W$ satisfies the {\em Calabi-Yau condition $\sum_i q_i=1$},
$X_W=\{W=0\}$ defines a Calabi-Yau hypersurface in the weighted projective space $\mathbb{P}^N(c_0,c_1, \dots, c_N)$, where $q_i=c_i/d$ for a common denominator $d$. The element $J_W$ acts trivially on $X_W$, while for any subgroup $G$ such that $\langle J_W\rangle \subseteq G\subseteq G_{max}$, the group $\widetilde{G}=G/\langle J\rangle$ acts faithfully on $X_W$.  The LG/CY correspondence predicts that the FJRW theory of $(W, G)$, up to analytic continuation and the quantization of a symplectic transformation, is equivalent to the Gromov-Witten theory of $X_{W}/\widetilde{G}$ \cite{R}. The case studied here are {\em mirror} of  the three classes of simple elliptic singularities: $\widetilde{E}_N(N=6,7,8)$.  More precisely,
\begin{align}
\notag
\P^1(3,3,3)=&\{P^T_8:=x^3_0+x^3_1+x^3_2=0\}/\widetilde{G}_{max},\\
\notag
\P^1(2,4,4)=&\{X^T_9:=x^2_0x_1+x^3_1+x_0x^2_2=0\}/\widetilde{G}_{max},\\
\notag
\P^1(2,3,6)=&\{J^T_{10}:=x^3_0+x^3_1+x_1x^2_2=0\}/\widetilde{G}_{max}.
\end{align}
Chiodo--Ruan (see \cite{CR1}), proposed a three-step approach to the LG/CY correspondence based on the B-model. Let us take simple elliptic singularities to simplify the notation. By Berglund-H\"ubsch-Krawitz, $P^T_8, X^T_9, J^T_{10}$ with $G_{max}$ is mirror to $P_8, X_9, J_{10}$ with the trivial group. There is a B-model construction in the latter case in terms of Saito-Givental theory. More precisely, consider the miniversal deformation of the simple elliptic singularities in the so-called {\em marginal direction}:
$P_8+s x_0x_1x_2, X_9+s x_0x_1x_2,$ or $J_{10}+s x_0x_1x_2$ for all the nonsingular values of $s$.
According to Saito, the above miniversal deformation space admits a generic semisimple Frobenius
manifold structure. Givental has constructed a higher genus generating function $\F_g$ over semisimple points.
We should mention that the original Saito-Givental theory is defined for a germ of singularities.
On the other hand, we study a "global" version of the Saito-Givental theory, where the marginal parameter $s$ is deformed from zero to infinity. In fact, the modularity arises only from this global point of view. To emphasis this key perspective,
we often refer to it as a {\em global Saito-Givental theory}.

Chiodo--Ruan (see \cite{CR1}) proposed that (i) FJRW theory of $(W, G_{max})$ is mirror to a
global Saito-Givental theory at $s=0$; (ii) GW theory of $X_W/\widetilde{G}_{max}$ is mirror to global Saito-Givental theory at $s=\infty$; (iii)
global Saito-Givental theory at $s=0$ is related to global Saito-Givental theory at $s=\infty$ by analytic continuation and quantization of a symplectic transformation. This article and that of Krawitz-Shen studied completely different aspects of this problem and can be treated as a single package. In particular, Krawitz and Shen proved (i), (ii) in \cite{KS} by a direct computation of the A-model for both GW theory and FJRW theory. In this article, we gave a detailed study of the B-model and affirm (iii). Therefore, by using our theorem (see Theorem \ref{desc:transf}), Krawitz-Shen deduced the LG/CY correspondence of all genera for the above examples.

Simple elliptic singularities are usually organized in one-parameter families, which according to Saito's interpretation (see \cite{Sa2})  can be viewed as a pull-back of some universal family parametrized by the modular curve.  Let us point out that we are slightly abusing notation, because for $X_9$ and $J_{10}$ we use respectively the normal forms $x_0^2x_2+x_0x_1^3+x_2^2$ and $x_0^3+x_1^3x_2+x_2^2$ instead of $x_0^4+x_1^4+x_2^2$ and $x_0^6+x_1^3+x_2^2$. The main motivation for our choice is to simplify the exposition. The LG/CY correspondence can be proved for other normal forms as well. We will deal with the remaining cases in a separate publication.

On the other hand, we proved much more than just (iii),  namely the modularity of global Saito-Givental theory! This is not a consequence of the LG/CY correspondence and represents an entirely new direction in GW theory. However, to draw the consequence for the A-model such as GW theory,  we use Krawitz-Shen's theorems at critical places. Namely, we use (i) to prove the extendibility of global Saito-Givental
theory over the caustic and (ii) to connect our result to GW theory.

Finally, the appearance of modularity in the B-model comes from the global behavior of the primitive form used to define the Frobenius structure. In general, it is a difficult problem to compute the primitive forms. However, in the cases under consideration, a primitive form and the corresponding Frobenius structure are determined by a choice of symplectic basis of $H_1$ of the corresponding elliptic curve. This basis determines a point $\tau\in \mathbb{H}$ on the upper half-plane for each value of the parameter $s$. The domain of the parameter $s$ can be identified with the quotient of $\mathbb{H}$ by the monodromy group $\Gamma$. Therefore, the global Saito-Givental function $\F_g$ should be viewed as a function of $\tau\in \mathbb{H}$. In this paper, we study the transformation of $\F_g(\tau)$ under $\tau\rightarrow \nu(\tau)$ for $\nu\in \Gamma$. The transformation of $\F_g(\tau)$ is given by the quantization of a certain symplectic transformation. This confirms (iii). We want to emphasise that $\F_g$ does \emph{not} transform as a modular form. A crucial idea, motivated by physics, is to complete $\F_g(\tau)$ in a specific way to a {\em non-holomorphic} function $\F_g(\tau, \bar{\tau})$. The anti-holomorphic completion is the generalization of the corresponding construction of quasi-modular form. Then, we can show
\begin{theorem}\label{t2}
The modified non-holomorphic function $\F_g(\tau, \bar{\tau})$ transforms as  an almost holomorphic modular form (see the detailed statement in Section \ref{quasi-modular}).
\end{theorem}
This implies that the original $\F_g(\tau)$ is quasi-modular. Using the results of Krawitz-Shen,  it implies Theorem \ref{t1}.

\subsection{Acknowledgements.} A special thanks goes to A. Klemm for sharing his insight and ideas on the subject of modularity and Gromov-Witten
theory. The second author would like to thank M. Aganagic and V. Bouchard for many interesting
 discussions. We are thankful to J. Lagarias for pointing out some very helpful references about modular forms and to I. Dolgachev and S. Galkin for many helpful discussions. We thank A. Greenspoon for editorial assistance. Both authors are supported by NSF grants.

%{\bf I eliminate the discussion here since we did not use later.}

%{\bf I suggest that we organize the paper in the follows: section two: Global Frobenius structure for simple elliptic singularities. Here, we discuss primitive form and calculate flat coordinates and so on. The central theorem is that all the data are analytic function of upper half plane $H^+$.

%section three: Givental function. Here, we introduce Givental's function and discuss its analyticity including at cusp.

 %section four: global Saito-Givental function and modular transformation. Here, we calculate formula under modular transformation.

 %section five: Anti-holomorphic modification and modularity. Here, we do our magic modification and prove the main theorem.

 %section six: Closed formula for $P_8$ case. We may or may not have this section depending on the future progress.

\section{ Global Frobenius structures}\label{global_fs}
To simplify the notation, we shall focus on  the simple elliptic singularities of the $P_8$-family:
\ben
f(\sigma,x)= x_0^3+x_1^3+x_2^3+\sigma x_0x_1x_2,\quad x=(x_0,x_1,x_2).
\een
$\sigma$ takes values in the punctured complex line $\Sigma=\{\sigma^3+27\neq 0\},$ so that the origin $x=0$ is an isolated critical point.
The remaining two cases can be analyzed in a similar way. The
necessary modifications are explained in Section \ref{PXJ}.

\subsection{Basic set-up}

   Let us first recall the basic set-up of Saito's Frobenius manifold structure on the miniversal deformation of a singularity. We will use our
   example to illustrate the procedure.

\subsubsection{The space of miniversal deformations}
Recall (see \cite{AGV}) the action of the group of germs of holomorphic changes of the coordinates $(\C^3,0)\to (\C^3,0)$ on the space of all germs at $0$ of holomorphic functions. Given a holomorphic germ $f(x)$ with an isolated critical point at $x=0$ we say that the family of functions $F(s,x)$ is a {\em miniversal deformation} of $f$ if it is transversal to the orbit of $f$. One way to construct a miniversal deformation is to choose a $\C$-linear basis $\{\phi_i(x)\}$ in the Jacobi algebra $\O_{\C^3,0}/\langle \d_{x_0}f,\d_{x_1}f,\d_{x_2}f\rangle$. Then the following family provides a miniversal deformation:
\ben
F(s,x)=f(x)+\sum_{i=1}^\mu s_i\phi_i(x),\quad s=(s_1,s_2,\dots,s_\mu)\in \C^{\mu},
\een
where $\mu$ is the dimension of the Jacobi algebra, also known as the {\em Milnor number} or the {\em multiplicity} of the critical point.

In our setting the Milnor number is $\mu=8$. It is convenient to use the index set $\{-1,0,1,\dots,6\}$ instead of $\{1,\dots,\mu\}.$ We choose the following set of monomials to construct a miniversal deformation:
$\phi_{-1}=x_0x_1x_2$, $\phi_0=1$, and $\phi_i$ $i=1,2,\dots,6$ are given respectively by:
$$
x_0,\quad x_1,  \quad  x_2,\quad x_0x_1,\quad x_0x_2,\quad  x_1x_2.
$$
\begin{comment}
\begin{align}\notag
\phi_{-1}(x)&=x_0x_1x_2,& \phi_0(x)&=1,& &\\
\notag
\phi_1(x)&=x_0,         & \phi_2(x)&=x_1,& \phi_3(x)&=x_2,\\
\notag
\phi_4(x)&=x_0x_1,      & \phi_5(x)&=x_0x_2,& \phi_6(x)&=x_1x_2.
\end{align}
\end{comment}
Note that $s_{-1}$ is naturally identified with $\sigma.$ Let us assign weight $1/3$ to each variable $x_i$ so that $f(\sigma,x)$ and $\phi_i(x)$ become weighted-homogeneous of degree respectively 1 and $1-d_i$, where
\ben
d_{-1}=0,\ d_0=1,\  d_1=d_2=d_3=2/3,\quad\mbox{and}\quad d_4=d_5=d_6=1/3.
\een
Note that assigning weight $d_i$ to each $s_i$ turns $F$ into a weighted-homogeneous function of weight 1.

Put $\S= \Sigma\times \C^{\mu-1}$ and $X=\S\times \C^3$. Then we have the following maps:
\ben
\begin{CD}
 \S\times \C^3 & & \\
@V{\varphi}VV\searrow &  \\
\S\times \C &@>>{p}> & \S
\end{CD}
\qquad
\begin{tabular}{rl}
&$\varphi(s,x)= (s,F(s,x))$,\\
&\\
&$p(s,\gl)=s.$
\end{tabular}
\een
By definition the critical set $C$ of $F$ is the support of the sheaf
\ben
\O_C:=\O_X/\langle \d_{x_0} F,\d_{x_1}F,\d_{x_{2}}F\rangle.
\een
The map $\d/\d s_i\mapsto \d F/\d s_i$ induces an isomorphism between the sheaf $\T_\S$ of holomorphic vector fields on $\S$ and $q_*\O_{C}$, where $q=p\circ\varphi$.
In particular, each tangent space $T_s\S$ is equipped with an associative commutative multiplication $\bullet_s$ depending holomorphically on $s\in \S$. If in addition we have a volume form $\omega=g(s,x)d^3x,$ where $d^3x=dx_0dx_1dx_2$ is the standard volume form; then $q_*\O_C$ (hence $\T_\S$ as well) is equipped with the {\em residue pairing}:
\beq\label{res:pairing}
(\psi_1,\psi_2) = \frac{1}{(2\pi i)^3} \ \int_{\Gamma_\ge} \frac{\psi_1(s,y)\psi_2(s,y)}{F_{y_0}F_{y_1}F_{y_2}}\, \omega,
\eeq
where $y=(y_0,y_1,y_2)$ are unimodular coordinates for the volume form, i.e., $\omega=d^3y$, and $\Gamma_\ge$ is a real 3-dimensional cycle supported on $|F_{x_0}|=|F_{x_1}|= |F_{x_2}|=\ge.$

Given a semi-infinite cycle
\beq\label{cycle}
\A\in  \lim_{ \longleftarrow } H_3(\C^3, (\C^3)_{-m} ;\C)\iso \C^\mu,
\eeq
where
\beq\label{level:lower}
(\C^3)_m=\{x\in \C^3\ |\ {\rm Re}(F(s,x)/z)\leq m\},
\eeq
put
\beq\label{osc_integral}
J_{\A}(s,z) = (-2\pi z)^{-3/2} \, zd_\S \, \int_{\A} e^{F(s,x)/z}\omega,
\eeq
where $d_\S$ is the de Rham differential on $\S$. The oscillatory integrals $J_\A$ are by definition sections of the cotangent sheaf $\T_\S^*$.

According to Saito's theory of primitive forms \cite{Sa1, SaT}, there exists a volume form $\omega$ such that the residue pairing is flat and the oscillatory integrals satisfy a system of differential equations, which in flat-homogeneous coordinates $t=(t_{-1},t_0,\dots,t_6)$ has the form
\beq\label{frob_eq1}
z\d_i J_\A(t,z) = \d_i \bullet_t J_\A(t,z),
\eeq
where $ \d_i:=\d/\d t_i\ (-1\leq i\leq 6)$ and the multiplication is defined by identifying vectors and covectors via the residue pairing. Due to homogeneity the integrals satisfy a differential equation with respect to the parameter $z\in \C^*$:
\beq
\label{frob_eq2}
(z\d_z + E)J_{\A}(t,z) =  \theta\, J_\A(t,z),
\eeq
where
\ben
E=\sum_{i=-1}^6 d_it_i \d_i,\quad (d_i:={\rm deg}\, t_i={\rm deg}\, s_i),
\een
is the {\em Euler vector field} and $\theta$ is the so-called {\em Hodge grading operator }:
\ben
\theta:\T^*_S\rightarrow \T^*_S,\quad \theta(dt_i)=\Big(\frac{1}{2}-d_i\Big)dt_i.
\een
The compatibility of the system \eqref{frob_eq1}--\eqref{frob_eq2} implies that the residue pairing, the multiplication, and the Euler vector field give rise to a {\em conformal Frobenius structure} of   conformal  dimension $1$. We refer to B. Dubrovin \cite{Du} for the definition and more details on Frobenius structures.

For the simple elliptic singularities of type $P_8$ the primitive forms can be described as follows. Let $\pi(\si)$ be a solution to the differential equation
\beq\label{Picard-Fuchs}
\frac{d^2u}{d\si^2} + \frac{3\si^2}{\si^3+27}\ \frac{du}{d\si} +  \frac{\si}{\si^3+27}\ u = 0;
\eeq
then the form $\omega=d^3x/\pi(\si)$ is primitive. For the reader's convenience we prove this statement in Appendix \ref{app:A}.

%In other words, if we have $\si_0\in \Si,$ such that $\pi(\si_0)\neq 0$, then we have a Frobenius structure defined on an open subset of $\S$ consisting of points $s$ such that $s_{-1}$ is sufficiently close to $\si_0.$

\subsection{Global Frobenius structures}

   Traditionally, one studies the germ at $s=0$ of the Frobenius structure in singularity theory, because in general the primitive form is known to exist only locally (as a germ with respect to the deformation parameters $s$). For our purposes however, we would like to vary the Frobenius structure from $s_{-1}=0$ to $s_{-1}=\infty$. This leads to the study of {\em global Frobenius structure.} In this subsection, we shall treat the construction of Saito's Frobenius manifold structure with this purpose in mind. Our first goal is to define primitive forms globally in the sense that they vary with $s$ in a nice fashion.

\subsubsection{Periods of elliptic curves and global primitive forms}

Put $X_{s,\gl}$ for the fiber $\varphi^{-1}(s,\gl)$ and let $(\S\times\C)'$ be the set of all $(s,\gl)$ that parametrize non-singular fibers $X_{s,\gl}$. The complement of $(\S\times\C)'$ is a analytic hypersurface, called {\em the discriminant}, and the union $X'$ of all non-singular fibers $X_{s,\gl}$ is a smooth fibration over $(\S\times\C)'$, called {\em the Milnor fibration}.  Following Looijenga \cite{Lo} we compactify the fibers $X_{s,\la}$ by adding an elliptic curve. Namely, the map
\ben
\S\times \C^3\to \S\times \C P^3,\quad (s,x)\mapsto (s,[x_0,x_1,x_2,1])
\een
is an embedding and we denote by $\overline{X}$ the Zariski closure of $X$ in $\S\times \C P^3$. The map $\varphi:X\to \S\times\C$ naturally extends to a map $\overline{X}\to \S\times\C$. We denote by $\overline{X}_{s,\la}$ the corresponding fibers. It is easy to check that the intersection of  $\overline{X}_{s,\la}$ with the hyperplane $\{z_3=0\}$ (here $[z_0,z_1,z_2,z_3]$ are homogeneous coordinates of $\C P^3$) is the elliptic curve (known also as the {\em elliptic curve at infinity}):
\ben
E_\sigma:\quad z_0^3+z_1^3+z_2^3+\sigma z_0z_1z_2 =0,
\een
where $\sigma=s_{-1}$. Moreover, the Gelfand--Leray form $d^3x/dF$ gives rise to a holomorphic form on $X_{s,\la}$ that has a simple pole along $E_\si$, and therefore its Poincar\'e residue ${\rm Res}_{E_\si}[d^3x/dF]$ is a holomorphic 1-form on $E_\si$ of degree 0, so it depends only  on $\si=s_{-1}$ but not on $s_0,s_1,\dots, s_6$ (see \cite{Lo}).

According to K. Saito (see \cite{Sa1}) the primitive forms for simple elliptic singularities are parametrized by the periods of $E_\sigma$
\beq\label{period}
\pi_A(\si):= 2\pi i \int_{A_\sigma} {\rm Res}_{E_\si}[d^3x/dF]\ ,
\eeq
where $A\in H_1(E_{\si_0},\C)$ is any non-zero 1-cycle and $A_\si$ is a flat family of cycles uniquely determined by $A$ for all $\si$ in a small neighborhood of $\si_0$. In Appendix \ref{app:A} we prove that the space of solutions to \eqref{Picard-Fuchs} coincides with the space of all periods $\pi_A(\si)$. Slightly abusing the notation, we often omit the index $\sigma$ from $A_\si$ and use $A$ to denote the flat family of cycles induced by $A$.

\begin{comment}
We make the following natural choice. Let $\si_0=0$ be the reference point. The fundamental group $\pi_1(\Si)$ acts naturally on the homology group $H_1(E_{\si_0},\Z)$ by monodromy transformations that form a group $\Gamma$. It is known  (see \cite{Ko}, \cite{Shi}) that $\Gamma$ is isomorphic to
\ben
\Gamma(3)=\Big\{
g\in {\rm SL}(2,\Z)\ |\ g\equiv I_2 ({\rm mod}\ 3)\Big\},
\een
where $I_n$ is the identity matrix of size $n$ and the congruence
relation is imposed on the corresponding entries of the matrices. In
particular, we can choose a cycle $A\in H_1(E_{\si_0};\Z)$ such that
it remains fixed  under the monodromy transformation that represents
the matrix $\begin{bmatrix}1& 3\\ 0& 1\end{bmatrix}$. Let us point out
that if we choose a different cycle $A$, not necessarily integral,
then the monodromy group $\Gamma$ becomes just a conjugate of
$\Gamma(3)$ in ${\rm SL}(2,\C)$ and the quasi-modularity property
(with respect to $\Gamma$) remains valid.
\end{comment}

For our purposes it is convenient to rewrite the integral (\ref{period}) as a period of the Gelfand--Leray form. Namely, let $X_{s,\la}$ be any non-singular fiber of the Milnor fibration such that $s_{-1}=\si$. The boundary of any tubular neighborhood of $E_\si$ in $\overline{X}_{s,\la}$ is a circle bundle over $E_\si$ that induces via pullback an injective {\em tube map} $L:H_1(E_\sigma)\to H_2(X_{s,\la})$. Let $a=L(A);$ then we have
\beq\label{lerey_period}
\pi_A(\si) = \pi_a(s):=\int_{a}\frac{d^3x}{dF}.
\eeq
We refer to $a$ as a {\em tube} or {\em toroidal} cycle. The space of
all toroidal cycles coincides with the kernel of the intersection
pairing on  $H_2(X_{s,\la};\C)$ (see \cite{Ga, Lo}).

\begin{comment}
In order to define the Frobenius structure globally, we have to choose
the cycle $A$ in such a way that $\pi_A(\si)\neq 0$ for all $\si\in
\Si$.  In some sense the main point of this paper is to describe the
 relation between the flat structures arising from two different
 choices of the cycle $A.$
\end{comment}

A flat family of cycles $A$ is a multi-valued object; therefore the induced global primitive form and global
Frobenius structure are multi-valued as well. This leads to the key
observation that, when discussing a global Frobenius structure, it is
more natural to replace $\Sigma$ by its universal cover. The latter is
naturally identified with the upper half-plane $\mathbb{H}.$ Namely,
fix a reference point, say $\si_0=0$. The points in the universal
cover $\tilde{\Sigma}$ of $\Sigma$ are pairs consisting of a point
$\sigma\in \Si$ and a homotopy class of paths $l(t)$ with $l(0)=\si_0, l(1)=\sigma.$ We fix a symplectic basis $\{A',B'\}$ of $H_1(E_{\si_0};\Z)$ once and for all. The map $(\si,l(t))\mapsto \tau'=\pi_{B'}/\pi_{A'}$, where the periods $\pi_{B'}$ and $\pi_{A'}$ are analytically continued along the path $l(t)$, defines an analytic isomorphism
between the universal cover of $\tilde{\Si}$ and the upper half-plane $\mathbb{H}$. In other words, we have a
Frobenius structure on $\mathbb{H}\times\C^{\mu-1}$ for any non-zero cycle
\beq\label{fstr:cycle}
A=dA'+cB'\in H_1(E_{\si_0};\C),\quad -d/c\notin\mathbb{H}.
\eeq

%%%%%%%%%%%%%%%%%%%%%%%%%%%%%%%%%%%%%%%%%%%%%%%%%%%

\subsubsection{Flat coordinates}
The goal in this section is to construct a flat homogeneous coordinate system $t=(t_{-1},t_0,\dots,t_6).$  The idea is to expand the oscillatory integrals into a power series near $z=\infty$. The flat coordinates will be identified with the leading coefficients in these expansions. The problem of analyzing the Gauss--Manin connection at $z=\infty$ was addressed by M. Noumi \cite{No} while the construction of flat coordinates for simple and simple elliptic singularities can be found in \cite{NY}. The combination of these two articles implies the result that we need. However, our point of view is somewhat different from the one in \cite{No}. For the reader's convenience as well as to avoid any misunderstanding we give a self-contained exposition.

Let $\ga(\si,1)\in H_2(X_{\si,1};\Z) $ be a flat family of cycles defined for $\si$ near $\si_0=0$. Using the rescaling $x\mapsto \la^{1/3}x$ we obtain a cycle $\ga(\si,\la) \in H_2(X_{\si, \la};\C)$. The cycle  $\A$ formed by $\ga(\si,z\la)$ as $\la$ varies along the path
\ben
\la:[0,\infty)\to \C,
\quad  \la(t)=-t,
\een
is a semi-infinite cycle of the type \eqref{cycle}. The corresponding oscillatory integral takes the form
\ben
\int_{\A} e^{F/z}\omega =z\, \int_0^{-\infty+i0} e^{\la} \ \int_{\ga(\si,\la)}e^{\sum_{j=0}^6 s_j\phi_j(x) z^{-d_j}}
\frac{\omega}{df} \ d\la.
\een
Rescaling $x\mapsto \la^{1/3}x$ and expanding the integrand in  powers of $z$ we get
\beq\label{expansion}
(-2\pi z)^{-3/2}\, \int_{\A} e^{F/z}\omega=z^{-1/2}\,\sum_{\delta} \Big(\int_{\ga(\si,1)} c_\delta(s,x)\frac{\omega}{df}\Big)\, z^{-\delta},
\eeq
where the sum is over all non-negative elements of the lattice in $\QQ$ spanned over $\Z$ by the degrees $d_i$ and  $c_\delta(s,x)$ is
\ben
\sum \widetilde{\Gamma}(k_0 (1-d_0)+k_1(1-d_1)+\cdots)\,
\frac{s_0^{k_0}}{k_0!} \frac{s_1^{k_1}}{k_1!}\cdots \, (\phi_0(x))^{k_0}(\phi_1(x))^{k_1}\cdots ,
\een
where the summation is over all integers  $k_i\geq 0$ such that  $k_0d_0+k_1d_1+\dots =\delta$ (note that the sum is finite) and
\ben
\widetilde{\Gamma}(k) :=(-2\pi)^{-3/2} \,(-1)^{k+1}\, \Gamma(k+1)=
(2\pi)^{-3/2} e^{\pi i(k-1/2)} \int_0^\infty e^{-t} t^k dt.
\een
Given a set of middle homology cycles $\ga_i$ ($-1\leq i\leq 6$) we define the matrix $\Pi$ whose entries are the following periods
\ben
\Pi_{\delta,i}=\int_{\ga_i(\si,1)} c_{\delta}(s,x)\frac{\omega}{df},\quad i=-1,0,\dots,6,
\een
where the index $\delta$ takes values in $\{0,1, 1/3,2/3\}$. The order in the latter set is such that it matches the rows in which the entries $\Pi_{\delta,i}$ should be placed.

Let $t=(t_{-1},t_0,\dots,t_6)$ be a flat-homogeneous coordinate system with degrees ${\rm deg} \, t_i=d_i$. It is convenient to introduce the following involution $'$ on the index set $\{-1,0,1,\dots,6\}$:
\beq\label{involution}
(-1)'= 0,\quad 0'=-1,\quad \mbox{and}\quad i'=7-i \quad \mbox{for}\quad 1\leq i\leq 6.
\eeq
Let us fix flat coordinates such that the residue pairing has the form $(\d_i,\d_j)=\delta_{ij'}$. Finally, we form the following matrix:
\beq\label{periods-matrix}
\begin{bmatrix}
t_{-1} & 1    \\
t^2/2  & t_0
\end{bmatrix} \oplus
\begin{bmatrix}
 0   &0   &0   &t_4 &t_5 & t_6 \\
 t_1 &t_2 &t_3 &0   &0   &0
\end{bmatrix},
\eeq
where $t^2=\sum_{i=-1}^6 t_it_{i'}.$ By direct sum $M_1\oplus M_2$ of two matrices $M_1$ and $M_2$ (not necessarily diagonal!) we mean a block-diagonal matrix with $M_1$ and $M_2$ on the diagonal.
\begin{lemma}\label{flat_coordinates}
There are cycles $\ga_i$ such that:
\begin{enumerate}
\item[(a)]
The period matrix $\Pi$ coincides with (\ref{periods-matrix}).
\item[(b)]
The cycles $\ga_i$ ($-1\leq i\leq 6$) are eigenvectors of the classical monodromy operator with eigenvalues
$e^{-2\pi\sqrt{-1}\,d_i}$.
\item[(c)]
The cycle $\ga_0=-(-2\pi)^{3/2}L(A)$, where $L(A)$ is the tube cycle that parametrizes the Frobenius structure.
\end{enumerate}
\end{lemma}
\proof
Let $\A_i$ be the semi-infinite cycles constructed from $\ga_i$ via rescaling. The oscillatory integrals $J_{\A_i}(s,z)$ satisfy the differential equations \eqref{frob_eq1} and \eqref{frob_eq2}.
The coordinates of $J_{\A_i}(s,z)$ with respect to the 1-forms $dt_{-1},dt_0,\dots,dt_6$ give rise to column vectors and we put  $J(s,z)$ for the matrix formed by these columns. Using \eqref{frob_eq2} we get that $J(s,z)$ has the following form:
\ben
(S_0 + S_1 z^{-1}+S_2z^{-2}+\cdots)\ z^{\theta},
\een
while \eqref{frob_eq1} implies that $S_0$ is a constant matrix independent of $t$ and $z$. Changing the cycles $\ga_i$ if necessary we can arrange that $S_0=1.$

\medskip

(a)
Let us compare the coefficients in front of $z^{-\delta}$ for $0\leq \delta\leq 1$ in
\beq\label{key_identity}
J_{\A_i}(t,z) = S(t,z)z^{\theta}dt_i.
\eeq
The RHS equals
\ben
z^{-d_{i}+1/2}\,S(t,z)dt_i=z^{1/2}\, (z^{-d_i}dt_i+\delta_{i,-1}z^{-1}S_1dt_{-1}+\cdots),
\een
where the dots stand for terms involving $z^{-\delta}$ with $\delta>1$. We have $S_1dt_{-1} = \sum t_i dt_{i'}$, because both co-vectors satisfy the differential equations $Lie_{\d_i} v(t) = dt_{i'}$ and the initial condition $v(0)=0$. Therefore, $S_1dt_{-1} = d t^2/2$. Comparing the coefficients in front of $z^{-\delta}$ for $0\leq \delta<1$ in \eqref{key_identity} we get (using also \eqref{expansion}) that $d_\S\Pi_{\delta,i}$ is either $0$ if $d_i\neq \delta$, or $dt_i$ if $d_i=\delta$. When $\delta=1$ we have:
\ben
d_\S\Pi_{1,0}=dt_0\quad \mbox{ and }\quad d_\S\Pi_{1,-1} = dt^2/2.
\een
In other words, up to some constant $4\times 8$ matrix $C$ the period matrix has the form that we want. In order to fix the constants we set $t_0=\dots=t_6=0$.
Up to some non-zero constant factors the differential forms $c_\delta(s,x)\omega/df$, $\delta=0,1,1/3, 2/3$, coincide repectively with
\ben
\frac{\omega}{df},\quad (s_0+\cdots)\frac{\omega}{df},\quad \Big(\sum_{i=4}^6s_i\phi_i(x)+\cdots\Big)\frac{\omega}{df},\quad
\Big(\sum_{i=1}^3s_i\phi_i(x)+\cdots\Big)\frac{\omega}{df},
\een
where the dots stand for at least quadratic polynomials in
$s_0, s_1,\dots, s_6.$ All periods vanish when $t_0=\dots=t_6=0$, except for $\Pi_{0,-1}$ and $\Pi_{0,0}$ (note that $\Pi_{0,i}=0,\ 1\leq i\leq 6$, follows from $\Pi_{1,i}=0$). We need to prove only that $C_{0,-1}=C_{0,0}=0.$ We return to these identities once we establish (b) and (c).

\medskip

(b) The Gelfand--Leray forms $\phi_i(x)\omega/df$ give rise to a basis of eigenvectors for the classical monodromy operator with eigenvalues $e^{2\pi \sqrt{-1}d_i}$. Since we already proved in (a) that $\Pi_{0,1}=\Pi_{1,1}=\Pi_{\frac{1}{3},1}=0$, we get
\ben
\int_{\ga_1}\phi_i\omega/df=0\quad  \mbox{for}\quad i=-1, 0, 4, 5, 6.
\een
In other words, $\ga_{1}$ belongs to the dual space of the space of middle cohomology classes spanned by $\phi_i(x)\omega/df$, $1\leq i\leq 3$. The latter is the eigenspace with eigenvalue $e^{2\pi\sqrt{-1}d_1}$; hence $\ga_1$ is an eigenvector with eigenvalue $e^{-2\pi\sqrt{-1}d_1}$. The remaining cases are analyzed in a similar way.

\medskip

(c)
 Let us substitute $t_1=\dots=t_6=0$ in the $2\times 2$ block of $\Pi$ formed by the intersection of the rows $d=0,1$ and the columns $i=-1,0$. We get the following table of identities:
\begin{align}\notag
&-(-2\pi)^{-3/2}\int_{\ga_{-1}(\si,1)}\frac{\omega}{df}=t_{-1}+C_{0,-1}& &
-(-2\pi)^{-3/2}\int_{\ga_{0}(\si,1)}\frac{\omega}{df}=1 + C_{0,0}\\
\notag & & & \\
\notag
& -(-2\pi)^{-3/2}s_0\int_{\ga_{-1}(\si,1)}\frac{\omega}{df}=t_0t_{-1} & &
-(-2\pi)^{-3/2}s_0\int_{\ga_{0}(\si,1)}\frac{\omega}{df}=t_0.
\end{align}
Put $\ga_0=m\,a+n\ga_{-1}$, where $a=L(A)$. Then the $(0,0)$-identity (keep in mind also the $(0,-1)$-identity) turns into
$$
-(-2\pi)^{-3/2}m+n(t_{-1}+C_{0,-1})=1+C_{0,0},
$$
It follows that $n=0$. The $(1,0)$-identity gives $-(-2\pi)^{-3/2}ms_0=t_0$, while the remaining two identities give $s_0(t_{-1}+C_{0,-1})=t_0t_{-1}=-(-2\pi)^{-3/2}ms_0t_{-1}.$ From here we get $m=-(-2\pi)^{3/2}$ and $C_{0,-1}=0$, which imply also that $C_{0,0}=0.$
\qed

\subsection{Modular transformations of the Frobenius structure. }\label{modular_transformations}

\begin{comment}
The primitive form and the corresponding Frobenius structure are single-valued only when $\sigma\in \Sigma$ varies in a small neighborhood of a fixed reference point $\si_0$.
\end{comment}

Every closed loop $C$ in $\Sigma$ based at $\si_0$ induces a monodromy transormation $\nu$ of both $H_2(X_{\si_0};\C)$ and $H^2(X_{\si_0};\C).$ We refer to $\nu$
as a {\em modular transformation}, while the set of all modular transformations forms a group which we call
the {\em modular group} of the family of singularities at hand.
Let us fix a basis of cycles $\ga_i$ satisfying the conditions in Lemma \ref{flat_coordinates}.
\subsubsection{Modular transformations}
The middle cohomology groups $H^2(X_{s,\la};\C)$ form a vector bundle equipped with a flat Gauss--Manin conncetion. Given a holomorphic form $\phi(s,x)d^3x$ the integrals $\int\phi(s,x)d^3x/dF$ and $\int d^{-1}(\phi(s,x)d^3x)$ determine naturally sections of the middle cohomology bundle. Here $d$ is the de Rham differential with respect to $x\in \C^3$ and $d^{-1}\omega$ means any $2$-form $\eta$ such that $d\eta=\omega.$ We have the following formulas for the covariant derivatives of such sections (see \cite{AGV}):
\ben
\nabla_{\d/\d s_i} \int d^{-1}(\phi(s,x)d^3x) = -\int \frac{\d F}{\d s_i}\,\phi(s,x)\,\frac{d^3x}{dF} + \int d^{-1}{\rm Lie}_{\d/\d s_i} (\phi(s,x)d^3x)
\een
and
\ben
\nabla_{\d/\d \la} \int d^{-1} (\phi(s,x)d^3x) = \int \phi(s,x)\frac{d^3x}{dF}.
\een
The second formula implies the following identity:
\beq\label{divergence}
\d_\la \int_\ga \phi(s,x)\d_{x_i}F \frac{d^3x}{dF} = \int_\ga \d_{x_i}\phi(s,x)\frac{d^3x}{dF},
\eeq
where $\ga$ is some middle homology cycle. Indeed, the integrand on the LHS equals $\phi(x)d^3x/dx_i$ while the one on the RHS is $d(\phi(x)d^3x/dx_i)/dF$.
\begin{lemma}\label{monodromy_N}
Let $\nu$ be a modular transformation; then the matrix of $\nu$ with respect to the basis $\{\ga_i\}_{i=-1}^6$ has the following block-diagonal form:
\beq\label{monodromy_matrix}
g\oplus{\rm Diag}(e^{2\pi i d_1k},\dots,e^{2\pi i d_6k}),
\eeq
for some $(g,k)\in {\rm SL}(2;\C)\times \Z.$
\end{lemma}
\proof
Let us compute the monodromy of the following sections $\int \phi_i(x)d^3x/df$, $-1\leq i\leq 6.$ The space spanned by the sections with $i=-1$ and $i=0$ is dual to the space of toroidal cycles, i.e., it is isomorphic to the homology group $H_1(E_\si; \C)$. Note that
\ben
\d_\si \int_\ga x_0\frac{d^3x}{df} =-\d_\la \int_\ga x_0^2x_1 x_2\frac{d^3x}{df}.
\een
On the other hand
\ben
x_0^2x_1x_2 = \frac{9}{\si^3+27}x_1x_2f_{x_0}+\frac{\si^2}{\si^3+27}x_0x_1f_{x_1}-\frac{3\si}{\si^3+27}x_1^2f_{x_2}
\een
Using formula (\ref{divergence}) we get (recall that $\phi_1(x)=x_0$)
\ben
\d_\si \int_\ga \phi_1(x)\frac{d^3x}{df} =
- \frac{\si^2}{\si^3+27}\int_\ga \phi_1(x)\frac{d^3x}{df}.
\een
Solving this differential equation for $\si$ we get
\ben
\int \phi_1(x)\frac{d^3x}{df} = (\si^3+27)^{-1/3}  A_1,
\een
where $A_1\in H^2(X_{\si_0,1};\C)$ is a flat section of the middle cohomology bundle. Under analytic continuation along a simple loop around $(-27)^{1/3}$ the RHS gains a factor of $\ge^{-1},$ where $\ge=e^{2\pi i/3}$. Therefore, $A_1$ is an eigenvector of the corresponding monodromy transformation with eigenvalue $\ge$. Similarly, one proves that
\ben
\int \phi_i(x)\frac{d^3x}{df} = (\si^3+27)^{-1+d_i}  A_i,\quad 1\leq i\leq 6,
\een
where $A_i$ ($1\leq i\leq 6$) are eigenvectors with eigenvalues $e^{-2\pi i d_i}.$

Finally, using Lemma \ref{flat_coordinates} we get
\ben
\int_{\ga_i}\phi_j(x)\frac{\omega}{df} = 0
\een
in the following two cases: (1) $i=1,2,3$ and $j=-1,0,4,5,6$;  (2) $i=4,5,6$ and
$j=-1,0,1,2,3$. This means that $\ga_i$ and $A_i$ belong to
eigenspaces that are dual to each other. The lemma follows.
\qed

\subsubsection{Modular transformations of the flat coordinates}
According to Lemma \ref{flat_coordinates} there is a cycle $B=bA'+aB'$, linearly independent from $A,$ s.t.,
\beq\label{modulus}
t_{-1}:= \frac{\pi_B}{\pi_A}= \frac{a\tau'+b}{c\tau'+d}
\eeq
is a flat coordinate and the residue pairing of the vector fields $1$ and $\d/\d t_{-1}$ is 1. More precisely, in the notation of Lemma \ref{flat_coordinates}, we have
\ben
\ga_{-1} = -(-2\pi)^{3/2}L(B),\quad \ga_0=-(-2\pi)^{3/2}L(A).
\een
The intersection pairing on the elliptic curve at infinity, up to a sign, is the same as the Seifert form of the corresponding toroidal cycles. Using Theorem 10.28(i) from \cite{He1}, we see that up to a sign the intersection number $A\circ B$ must be $\sqrt{-1}$.
\begin{remark} The basis $\{A,B\}$ is not symplectic and $t_{-1}$ is not a modulus of the elliptic curve at
infinity.
\end{remark}
The analytic continuation along a path $C$ transforms $t_{-1}$ into
\beq\label{action_on_tau}
g(t_{-1}):=\frac{n_{11}t_{-1}+n_{21}}{n_{12}t_{-1}+n_{22}},
\eeq
where $(n_{ij})$ is the matrix (from $\ {\rm SL}_2(\C)$) that describes the parallel transport $g$ of
$\{\ga_{-1},\ga_0\}$ along $C$, i.e.,
\beq\label{nu}
g(\ga_{-1})=n_{11}\ga_{-1}+n_{21}\ga_0,\quad\mbox{and}\quad
g(\ga_{0})=n_{12}\ga_{-1}+n_{22}\ga_0.
\eeq
For a given matrix $g=(n_{ij})\in {\rm SL}_2(\C)$ we adopt the number theorist's notation
$$
j(g,t_{-1}):=n_{12}t_{-1}+n_{22}.
$$
\begin{lemma}\label{monodromy_transformation}
The analytic continuation along the loop $C$ induces a coordinate change $t\mapsto \nu(t)$ of the following form:
\ben
\nu(t)_{-1}= g(t_{-1}),\quad \nu(t)_{0}= t_0+\frac{n_{12}}{2j(g,t_{-1})}\ \sum_{i=1}^6 t_it_{i'},\quad \nu(t)_{i}=\frac{e^{2\pi i d_ik}}{j(g,t_{-1})}\ t_i \ (1\leq i\leq 6),
\een
where $k$ is some integer.
\end{lemma}
\proof
Note that the cycles $\ga_{-1}$ and $\ga_0$ are transformed respectively into $n_{11}\ga_{-1}+n_{21}\ga_0$ and $n_{12}\ga_{-1}+n_{22}\ga_0$. The period $\int_{\ga_0}d^3x/df$ is transformed into $(n_{12}t_{-1} + n_{22})\int_{\ga_0}d^3x/df$, which implies that the primitive form $\omega$ transforms into $\omega/j(\nu,t_{-1})$. According to Lemma  \ref{flat_coordinates}, we have $t_0=\int_{\ga_0}c_1(s,x)\omega/df$. Hence $t_0$ is transformed into
\ben
\Big( n_{22}t_0 + n_{12}\int_{\ga_{-1}}c_1(s,x)\omega/df\Big) j(\nu,t_{-1})^{-1}.
\een
According to Lemma \ref{flat_coordinates} the above integral is
$\Pi_{1,-1}=t_{-1} t_0 +\frac{1}{2} \sum_{j=1}^6 t_j t_{j'}$. This
proves the transformation law for $t_0$. The remaining ones are proved
in a similar way by using Lemmas \ref{flat_coordinates} and
\ref{monodromy_N}. \qed

\subsection{Changing the Frobenius structures}\label{cover}
Let $A_i,$ $i=1,2$ be two cycles of the form \eqref{fstr:cycle}. We pick a cycle $B_i$ for $A_i$ as explained above
and let $g_i$ be the matrix such that $t_{-1}^{i}=g_i(\tau').$ Since the intersection numbers $A_i\circ B_i(i=1,2)$
are equal, the matrix
\ben
g:=g_2g_1^{-1}=:\begin{bmatrix}n_{11}&n_{21}\\n_{12}&n_{22}\end{bmatrix}\in {\rm SL}_2(\C).
\een
Let $k\in \Z$ be an arbitrary integer.
\begin{lemma}\label{flat:coord}
The map  $t\mapsto \widetilde{t}$, defined by
\ben
\widetilde{t}_{-1}=g(t_{-1}),\quad \widetilde{t}_0= t_0+\frac{n_{12}}{2j(g,t_{-1})}\sum_i t_it_{i'},\quad
\widetilde{t}_i=\frac{e^{2\pi i d_ik}}{j(g,t_{-1})}\, t_i,
\een
respects the residue pairings corresponding to $A_1$ and $A_2$.
\end{lemma}
The proof is straightforward and it is omitted. In particular, this lemma allows us to identify the flat
vector fields arising from two different families of flat cycles $A_1$ and $A_2$. Note however, that
the corresponding Frobenius multiplications are identical if and only if the cycle $A_2$ is obtained from $A_1$ by means of parallel transport along a closed loop. The reason for this is that every analytic isomorphism $t_{-1}^{1}=g_1(\tau')\mapsto t_{-1}^{2}=g_2(\tau')$ has the form $g_2gg_1^{-1}$, where $g$  is an automorphism of
$\mathbb{H}$, i.e., $g\in {\rm SL}_2(\mathbb{R}  ).$ The structure constants of the Frobenius multiplications are functions of $\si$; therefore if we think of $\si$ as a function on the upper half-plane, then it should be $g$-invariant.
But the automorphisms of $\mathbb{H}$ that preserve $\si$ are precisely the elements of the modular group of the
$P_8$-singularity, i.e., the group of deck transformations of the universal cover $\mathbb{H}\to \Si.$

\section{Givental's higher genus potential}

In this section, we introduce Givental's higher genus potential $\F_{g, formal}$ which will be our object of study.

\subsection{Symplectic vector space}

Let $H$ be the space of flat vector fields on $\mathcal{S}$ equipped with the residue pairing $(\ ,\ )$. Following Givental we introduce the vector space $\H=H((z))$ of formal Laurent series in $z^{-1}$ equipped with the symplectic structure
\ben
\Omega(f(z),g(z))= {\rm res}_{z=0} (f(-z),g(z))dz.
\een
Using the polarization $\H=\H_+\oplus \H_-$, where $\H_+=H[z]$ and $\H_-=H[[z^{-1}]]z^{-1}$ we identify $\H$ with the cotangent bundle $T^*\H_+$. Let us fix flat coordinates
$$
t=(t_{-1},t_0,\dots,t_6),\quad  (\d_i,\d_j)=\delta_{i,j'},
$$
where $\d_i=\d/\d t^i$ and $'$ is the involution (\ref{involution}).

Using the residue pairing we identify the tangent  and the cotangent
bundle $T\mathcal{S}\cong T^*\mathcal{S}$.
%The connection defined by (\ref{frob_eq1}) and (\ref{frob_eq2}) induces a flat connection on $T^*\mathcal{S}$, i.e., the following system of differential equations is compatible:
%\beqa
%\label{eq1}
%z\d_i J(t,z) & = & \d_i\bullet_t J (t,z),\quad -1\leq i\leq 6, \\
%\label{eq2}
%(z\d_z+E) J (t,z) & = & \theta J(t,z).
%\eeqa
Using the flat coordinates we trivialize the cotangent bundle
$T^*\mathcal{S}\cong \mathcal{S}\times H$. In this way, $H$ turns into
the space of flat holomorphic differential 1-forms. We use the basis
$\{dt_i\}_{i=-1}^6$ of $H$ in order to represent the linear transformations of $H$
by matrices of size $\mu.$

\subsubsection{The stationary phase asymptotics}\label{asymptotic:sec}

Let $s\in \S$ be a semi-simple point, i.e., the critical values $u_i$ ($1\leq i\leq \mu$) form locally near $s$ a coordinate system. Then we have an isomorphism
\ben
\Psi:\C^\mu\to H\cong T_s\S,\quad e_i\mapsto \sqrt{\Delta_i}\,\d/\d u_i,\quad (\d/\d u_i,\d/\d u_j)=\delta_{ij}/\Delta_i,
\een
that diagonalizes the Frobenius multiplication and the residue pairing:
\ben
e_i\bullet e_j = \sqrt{\Delta_i}e_i\delta_{i,j},\quad (e_i,e_j)=\delta_{ij}.
\een
The system of differential equations (\ref{frob_eq1}) and (\ref{frob_eq2}) admits a unique formal solution of the type
\ben
\Psi R(s,z) e^{U/z},\quad R(s,z)=1+R_1(s)z+R_2(s) z^2+\cdots
\een
where $U$ is a diagonal matrix with entries $u_1,\dots,u_\mu$ on the
diagonal and $R_k(s)\in {\rm End}(\C^\mu)$.
Alternatively this formal
solution coincides with the stationary phase asymptotics of the
following integrals. Let $\B_i$ be the semi-infinite cycle of the
type \eqref{cycle} consisting of all points $x\in \C^3$ such that the gradient trajectories of $-{\rm Re}(F/z)$ flow into the critical value $u_i$. Then
\ben
(-2\pi z)^{-3/2}\ zd_S\ \int_{\B_i} e^{F(s,x)/z}\omega \ \sim\ e^{u_i/z}\Psi R(s,z) e_i\quad \mbox{ as }\ z\to 0.
\een
We refer to \cite{AGV, G1} for more details and proofs.

\subsection{The total ancestor potential}

Let us fix the Darboux coordinate system on $\H$ given by the linear functions $q_k^i$, $p_{k,i}$ defined as follows:
\ben
\f(z) = \sum_{k=0}^\infty \sum_{i=-1}^6\ (q_k^i\, \d_i\, z^k + p_{k,i}\,dt_i\,(-z)^{-k-1})\quad \in \quad \H,
\een
where $dt_i$ is identified via the residue pairing with $\d_{i'}$.

It is known (and it is easy to prove) that $R$ is a symplectic transformation, i.e., $\leftexp{T}R(-z)R(z)=I_\mu.$ Note that ${R}$ has the form $e^{A(z)},$ where $A(z)$ is an infinitesimal symplectic transformation. On the other hand, a linear transformation $A(z)$ is infinitesimal symplectic if and only if the map $\f\in \H \mapsto A\f\in \H$ defines a Hamiltonian vector field with Hamiltonian given by the quadratic function $h_A(\f) = \frac{1}{2}\Omega(A\f,\f)$. By definition, the quantization of $e^A$ is given by the differential operator $e^{\widehat{h}_A},$ where the quadratic Hamiltonians are quantized according to the following rules:
\ben
(p_{k,i}p_{l,j})\sphat = \hbar\frac{\d^2}{\d q_k^i\d q_l^j},\quad
(p_{k,i}q_l^j)\sphat = (q_l^jp_{k,i})\sphat = q_l^j\frac{\d}{\d q_k^i},\quad
(q_k^iq_l^j)\sphat =q_k^iq_l^j/\hbar.
\een
Note that the quantization defines a projective representation of the Poisson Lie algebra of quadratic Hamiltonians:
\ben
[\widehat{F},\widehat{G}]=\{F,G\}\sphat + C(F,G),
\een
where $F$ and $G$ are quadratic Hamiltonians and the values of the cocycle $C$ on a pair of Darboux monomials is non-zero only in the following cases:
\beq\label{cocycle}
C(p_{k,i}p_{l,j},q_k^iq_l^j)=
\begin{cases}
1 & \mbox{ if } (k,i)\neq (l,j),\\
2 & \mbox{ if } (k,i)=(l,j).
\end{cases}
\eeq
The action of the operator $\widehat{R}$ on an element $F(\q)\in \C_\hbar[[q_0,q_1+1,q_2,\dots]]$, whenever it makes sense, is given by the following formula:
\beq\label{R:fock}
\widehat{R}\,F(\q) = \left.\Big(e^{\frac{\hbar}{2}V\d^2}F(\q)\Big)\right|_{\q\mapsto R^{-1}\q},
\eeq
where $V\d^2$ is the quadratic differential operator $\sum_{k,l}(\d^a,V_{kl}\d^b)\d_{ q_k^a}\d_{ q_l^b}$, whose coefficients $V_{kl}$ are given by
\beq\label{V:fock}
\sum_{k,l=0}^\infty V_{kl}(-z)^k(-w)^l = \frac{\leftexp{T}{R}(z)R(w)-1}{z+w}.
\eeq

By definition, the Kontsevich--Witten tau-function is the following generating series:
\beq\label{D:pt}
\D_{\rm pt}(\hbar;q(z))=\exp\Big( \sum_{g,n}\frac{1}{n!}\hbar^{g-1}\int_{\overline{\M}_{g,n}}\prod_{i=1}^n (q(\psi_i)+\psi_i)\Big),
\eeq
where $q(z)=\sum_k q_k z^k,$ $(q_0,q_1,\ldots)$ are formal variables, $\psi_i$ ($1\leq i\leq n$) are the first Chern classes of the cotangent line bundles on $\overline{\M}_{g,n}.$ The function is interpreted as a formal series in $q_0,q_1+1,q_2,\ldots$ whose coefficients are Laurent series in $\hbar$..

Let $s\in \mathcal{S}$ be {\em a semi-simple} point, i.e., the critical values $u^i(s)$ ($1\leq i\leq \mu$) of $F(s,x)$ form a coordinate system. Let $t=(t_{-1},t_0,\dots,t_6)$ be the flat coordinates of $s$. Motivated by the Gromov--Witten theory of symplectic manifolds Givental introduced the notion of the {\em total ancestor potential} of a semi-simple Frobenius structure. In particular the definition makes sense in singularity theory as well. Namely, the total ancestor potential is by definition the following formal function on $H[z]$:
\beq\label{ancestor}
\A_t(\hbar;\q) :=
\widehat{\Psi}\,\widehat{R}\,e^{\widehat{U/z}}\,
\prod_{i=1}^\mu \D_{\rm pt}(\hbar\,\Delta_i; \leftexp{\it i}{\bf q} (z) \sqrt{\Delta_i})
\eeq
where
$$
\q(z)=\sum_{k=0}^\infty\sum_{a=-1}^6 q_k^az^k\d_a, \quad
  \leftexp{\it i}{\bf q} (z) = \sum_{k=0}^\infty \leftexp{\it i}{q}_k z^k.
$$
The quantization $\widehat{\Psi}$ is interpreted as the change of variables
\beq\label{change}
\sum_{i=1}^\mu \leftexp{\it i}{\bf q}(z)e_i=\Psi^{-1}\q(z)\quad \mbox{i.e.}\quad
\leftexp{\it i}{q}_k\sqrt{\Delta_i} =\sum_{a=-1}^6 (\d_a u^i)\, q_k^a.
\eeq
The correctness of definition \eqref{ancestor} is not quite obvious. The problem is that the substitution $\q\mapsto R^{-1}\q $, which,  written in more detail, reads
\ben
q_0\mapsto q_0,\quad q_1\mapsto \overline{R}_1q_0 + q_1,\quad q_2\mapsto \overline{R}_2q_0 +\overline{R}_1q_1 + q_2,\quad \dots,
\een
where
$$
R^{-1}=1+\overline{R}_1z+\overline{R}_2 z^2+\cdots,
$$
is not a well-defined operation on the space of formal series.
This complication however, is offset by a certain property of the Kontsevich--Witten tau function. By definition, an {\em asymptotic} function is an element of the Fock space $\C_\hbar[[q_0,q_1,\dots]]$ of the form
\ben
\A=\exp \Big(\, \sum_{g=0}^\infty \overline{\F}^{(g)}(\q)\hbar^{g-1}\,\Big).
\een
It is called {\em tame} if the following $(3g-3+r)$-jet constraints are satisfied:
\ben
\left. \frac{\d^{r}\F^{(g)}}{\d q_{k_1}^{i_1}\cdots \d q_{k_r}^{i_r}}\right|_{\q = 0} = 0 \quad \mbox{ if } k_1+\cdots + k_r\geq 3g-3+r.
\een
The Kontsevich--Witten tau function, up to the shift $q_1\mapsto q_1+{\bf 1}$, is tame for dimensional reasons: ${\rm dim}\, \mathcal{M}_{g,r}=3g-3+r$. It follows that the action of $\widehat{R}$ is well-defined. Moreover, according to Givental \cite{G2}, $\widehat{R}$ preserves the class of tame asymptotic functions. In other words, the total ancestor potential is a tame asymptotic function in the Fock space $\C_\hbar[[q_0,q_1+{\bf 1},q_2,\dots]].$

 Let us point out the following homogeneity property of the Kontsevich--Witten tau-function:
\ben
\D_{\rm pt}(c^2\hbar ;cQ(z)) = c^{-1/24}\ \D_{\rm pt}(\hbar;Q(z))\quad \mbox{ for all } c\in \C.
\een
This formula can be used to rewrite the definition of the ancestor
potential \eqref{ancestor} in a different form, which looks
simpler. However, we prefer to work with formula \eqref{ancestor}, otherwise
$\A_t$ will be a formal series in a different Fock space.

\subsubsection{The total ancestor potential at non-semi-simple points}
Equation \eqref{frob_eq2} can be rewritten as
\ben
\nabla_t J=0,\quad \mbox{where}\quad \nabla_t:= d - z^{-1}\theta + z^{-2}E\bullet_t\ .
\een
One may think of $\nabla_t$ as an {\em isomonodromic} family of connection operators on $\C\setminus{\{0\}}$ parametrized by $t\in \mathcal{S}$. Let $S(t,z)$ be gauge transformations of the form
\ben
1+ S_1(t)z^{-1}+S_2(t)z^{-2}+\cdots,
\een
conjugating $\nabla_t$ and $\nabla_0=d-z^{-1}\theta$:
\ben
\nabla_t=S(t,z)\,\nabla_0\, S(t,z)^{-1}.
\een
The series $S(t,z)$ is also a symplectic transformation so it can be quantized in the same way as $R$.  The quantized symplectic transformation $\widehat{S}$ acts as follows:
\beq\label{S:fock}
\widehat{S}^{-1}\, F(q) = e^{W(\q,\q)/2\hbar}F([S\q]_+),
\eeq
where $W(\q,\q)$ is the quadratic form $\sum_{k,l}(W_{kl}q_l,q_k)$ whose coefficients are defined by
\beq\label{W:fock}
\sum_{k,l\geq 0} W_{kl}z^{-k}w^{-l}=\frac{\leftexp{T}{S}(z)S(w)-1}{z^{-1}+w^{-1}}.
\eeq
The $+$ sign in \eqref{S:fock} means truncation of all negative powers of $z$, i.e., in $F(\q)$ we have to substitute:
\ben
q_k\mapsto q_k + S_1q_{k+1}+S_2q_{k+2}+\cdots,\quad k=0,1,2,\dots\ .
\een
This operation is well-defined on the space of formal series. Note however that $S_1{\d_0}=t$ (see the proof of Lemma \ref{flat_coordinates}), where $t=\sum t_a\d_a\in H$ are the flat coordinates of the point $s\in \mathcal{S}$. Therefore, we have an isomorphism
\ben
\widehat{S}^{-1}:
\C_\hbar[[q_0,q_1+{\bf 1},q_2,\dots]]\to \C_\hbar[[q_0-t,q_1+{\bf 1},q_2,\dots]].
\een
Following Givental, we define the so-called {\em total descendant potential}:
\beq\label{DAN}
\D(\hbar;\q)= e^{F^{(1)}(t)}\, \widehat{S(t,z)}^{-1}\,\A_t(\hbar;\q),
\eeq
where
\beq\label{genus1}
F^{(1)}(t):= \frac{1}{2}\sum_{i=1}^\mu \int R_1^{ii}du^i + \frac{1}{48}\sum_{i=1}^\mu \ \ln(\Delta_i),
\eeq
is called the {\em genus}-1 potential. Since $S(t,z)$ and $\Psi R e^{U/z}$ satisfy the same differential equations with respect to $t$, one can check that the definition  \eqref{DAN} is independent of $t$, i.e., $\d_i\D=0$ for all $i$. By setting $t=q_0$, we get that the total descendant potential is a formal series in $q_1+{\bf 1},q_2,q_3,\dots$, whose coefficients are analytic  multi-valued functions on $\S$ with poles along the {\em caustic} $\mathcal{K}$. Here, multi-valued means that they are single-valued on the universal cover of $\S$, while the caustic is the subset of $\S$ of all non-semi-simple points.

Since the calibration $S$ is defined for all $s\in \S$, we can use equation
 \eqref{DAN} to define $\A_t$ for all $t$ as well. Note however that in the setting of an arbitrary semi-simple Frobenius structure $\A_t$ might
 not be a power series in $q_0$. In the setting of singularity theory, Givental (see \cite{G1}) conjectured that
\begin{conjecture}\label{hol:ext} The coefficients of the total descendant potential $\D(\hbar;\q)$ extend holomorphically through the caustic $\mathcal{K}.$
\end{conjecture}
In particular, if this is true then the total ancestor potential $\A_t$ is a power series in $q_0,q_1+{\bf 1},q_2,\dots$ whose coefficients are holomorphic functions in $t$.

\subsection{Extending $\A_t$ over the non-semi-simple locus}
The primitive form is multi-valued on $\Sigma$, but it is  analytic on the universal cover
$\widetilde{\Sigma}\cong \mathbb{H}$ of $\Sigma$ (see Subsection \ref{cover}). Therefore,
the Frobenius structure, which a priori is defined only for $s\in \mathcal{S}$ such that $s_{-1}$
is near $\si_0$, induces a holomorphic Frobenius structure on $\mathbb{H}\times\C^{\mu-1}$.
Let $\widetilde{\mathcal{K}}$ be the lift of the caustic to the universal cover, i.e., the set of all $t\in\mathbb{H}\times\C^{\mu-1}$ such that    the critical values $\{u^i(t)\}_{i={-1}}^6$ fail to form a coordinate system.  Then the total ancestor potential is a formal series whose coefficients are holomorphic on $\mathbb{H}\times\C^{\mu-1}\setminus \widetilde{\mathcal{K}}.$
\begin{lemma}\label{analyticity} Assume that the coefficients of the
  total ancestor potential are holomorphic in a neighborhood of some
  point $\tau_0\times 0\in \mathbb{H}\times\C^{\mu-1}$. Then they extend holomorphically across the caustic $\widetilde{\mathcal{K}}$.
\end{lemma}
\proof
Let $a(t)$ be one of the coefficients of $\A_t.$ Since the operators
$R_k$ and the Hessians $\Delta_i$ have only finite order poles along
$\widetilde{\mathcal{K}}$ the same is true  for $a(t)$. In other words the set
$\widetilde{\mathcal{K}}_a$ of all points $t\in \mathbb{H}\times \C^{\mu-1}$ such
that $a(t)$ is not holomorphic is a analytic subset. Let us assume
that $\widetilde{\mathcal{K}}_a$ is non-empty. Due to the Hartogs extension
theorem the codimension of  $\widetilde{\mathcal{K}}_a$ is at least 1 and hence it
is exactly 1. According to the assumption of the lemma
$\mathbb{H}\times {0}$ is not contained in $\widetilde{\mathcal{K}}_a$. It follows
that $\widetilde{\mathcal{K}}_a$ intersects  $\mathbb{H}\times {0}$ in a discrete
subset $\{\tau_i\times 0\}$. Moreover, due to homogeneity
$\widetilde{\mathcal{K}}_a$ is invariant with respect to the rescaling action
(with approriate weights) of $\C^*$ on $\mathbb{H}\times\C^{\mu-1}.$
Therefore every irreducible component of $\widetilde{\mathcal{K}}_a$ intersects
$\mathbb{H}\times {0}$, because the coordinates on $\C^{\mu-1}$ have
positive weights, so every $\C^*$-orbit intersects $\mathbb{H}\times
{0}$. It follows that $\widetilde{\mathcal{K}}_a$ is a disjoint union of
irreducible components of the type $\{\tau_i\}\times \C^{\mu-1}$. In
particular, the caustic $\widetilde{\mathcal{K}}$ has irreducible components of
this type as well. But this is not true: it is easy to see that
$\{\tau_i\}\times \C^{\mu-1}$ has a semi-simple point for every
$\tau_i$.
\qed

\medskip

According to M. Krawitz and Y. Shen \cite{KS}, if we choose the cycle $A$ in
such a way that $\pi_A(0)=1$ and $\pi_A'(0)=0$ then the total
descendant potential coincides with the generating function for the
FJRW invariants of the singularity
$f(0,x)=x_0^3+x_1^3+x_2^3$. Then the assumptions of the lemma are
satisfied with $\si_0=0$. Now if we choose
a different primitive form, then it is easy to see that the total ancestor potential
changes by a formula similar to the one in Theorem \ref{desc:transf}.  In particular,
the new potential depends holomorphically on $t$ as well.

\subsection{The poles at the cusps}

To close the section, we prove that
\begin{lemma}
The coefficients of the
  total ancestor potential have at most finite order poles at the cusps.
  \end{lemma}
  \proof
Recall the notation from Section \ref{asymptotic:sec}. The coefficients $R_k$ are determined recursively by the following relations:
\ben
(d+\Psi^{-1}d\Psi\wedge)R_k=[dU,R_{k+1}],
\een
which determines the off-diagonal entries of $R_{k+1}^{ij}$ in terms of the entries of $R_k$, and
\ben
R_{k+1}^{ii}=\frac{1}{k+1}\ \sum_{j\neq i} \ R_1^{ij}\,R_{k+1}^{ji}\ (u_i-u_j).
\een
These formulas are derived from the fact that the asymptotic operator $\Psi R e^{U/z}$ is a solution to the system of differential equations \eqref{frob_eq1} and \eqref{frob_eq2} (see \cite{G1}). In order to prove that the coefficients of the ancestor potential have finite order poles at the cusps, it is enough to prove that the asymptotic operator has finite order poles at $\si=\infty$ and $\si^3+27=0.$

Given a point $t=(t_{-1},t_0,\dots,t_{\mu-2})$ we put  $'t=(t_0,\dots,t_{\mu-2})$ and view the critical values as functions in $(\si,'t)\in \Si\times \C^{\mu-1}$. We need to prove that (for $'t$ fixed) $u_i(\si,'t)$ has a finite order pole at the punctures of the Riemann sphere $\Si$. Let us show how the argument works for one of the finite punctures $\si_0$, i.e.,  $\si_0$ is such that $\si_0^3+27=0$. For the puncture at $\si=\infty$ the argument is similar.

It is well known that the critical values are eigenvalues of the multiplication by the Euler vector field, i.e., they are the zeroes of an algebraic equation
\ben
u^\mu+\sum_{k=1}^\mu a_k(\si,\ 't)u^{\mu-k} :={\rm det}(u\, I_\mu - E\bullet_t )=0.
\een
It is easy to see that there is an integer $m$ and a constant $C$ such that
\ben
|a_k(\si,\ 't)|<C\, (\si^3+27)^{-m}
\een
for all $(\si,\ 't)$ in a fixed neighborhood of $(\si_0,0)$. Since we have
\ben
1+\sum_{k=1}^\mu a_k(\si,\, 't) u_i(\si,\ 't)^{-k} =0
\een
at least for one $k$ we must have
\ben
|a_k(\si,\, 't) u_i(\si,\ 't)^{-k}|\geq 1/\mu.
\een
From here one gets easily that
\ben
|u_i(\si,\ 't)\ (\si^3+27)^{m/k}|\leq C\mu.\qed
\een

\section{Transformations of the ancestor potentials}
\label{quasi-modular}
Let us fix a flat coordinate system
$$
t=(t_{-1},t_0,\dots,t_6),\quad (\d_i,\d_j)=\delta_{ij'},
$$
corresponding to an arbitrary primitive form.  It is convenient to denote the remaining coordinates by $'{t}=(t_0,t_1,\dots,t_6)$. Slightly abusing the notation we will sometimes identify $t_{-1}$ with the point $(t_{-1},0)$.

Recall that analytic continuation along some closed loop in $\Si$ transforms the flat coordinates $t\mapsto \nu(t)$ according to the formulas in Lemma \ref{monodromy_transformation}. In this section we would like to calculate how the total ancestor potential $\A_{t_{-1}}(\hbar;\q)$ transforms under analytic continuation.
 We will assume that Conjecture \ref{hol:ext} holds, so that we can view the total ancestor potential as a formal series in $q_0,q_1+{\bf 1},q_2,\dots$. This assumption is not really necessary in order to prove the transformation law of  $\A_{t_{-1}}$, but it is necessary later on in order to prove that the coefficients of $\A_{t_{-1}}$ are quasi-modular forms.

\subsection{Modular transformations}
We start by determining how the operator $\Psi R e^{U/z}$ changes under analytic continuation along some closed loop in $\Si$. Let $\nu$ be the corresponding modular transformation of the middle homology group $H_2(X_{\si_0,1};\Z)$. Note that if we fix a Morse coordinate system near each critical point $\xi_i$, then analytic continuation will simply permute the basis. Hence, the monodromy transformation of the stationary phase asymptotics of the oscillatory integrals $\int_{\B_i}e^{F/z}\omega$ is represented by some permutation matrix $P$.
Finally, given $\nu=(g,k)\in {\rm SL}_2(\C)\times\Z$, put
\beq\label{matrix:M}
M_{\nu}=\begin{bmatrix}
j(g,t_{-1})^{-1}  & *               &   *        &  *   \\
 0               & j(g,t_{-1})     &   0        &  0   \\
 0               & *               & \ge^{2k} I_3  &  0    \\
 0               & *               &   0        & \ge^k I_3
\end{bmatrix}
\eeq
where the $*$ entries in the first row and the second column are respectively
\ben
M_{-1,j}= -e^{2\pi i d_jk}\,n_{12}j(g,t_{-1})^{-1}t_{j},\quad 1\leq j\leq 6
\een
and
\ben
M_{-1,0}= -n_{12}z - \frac{n_{12}^2}{2j(g,t_{-1})}\sum_{i=1}^6 t_it_{i'},\quad  M_{i,0}=n_{12}t_{i'}  \quad 1\leq i\leq 6,
\een
$\ge=e^{2\pi i/3}.$

Let us point out that the transposition of a given matrix $A$ with respect to the residue pairing has
the following form:
\beq\label{transposition}
(\leftexp{T}{A})_{ij}=A_{j'i'}, \quad -1\leq i,j\leq 6.
\eeq
\begin{lemma}\label{ancestor:mon}
Analytic continuation changes the operator $\Psi R e^{U/z}$ into
$$
 \leftexp{T}{M}_\nu\,(\Psi R e^{U/z})\, P,
$$
where $\nu=(g,k)$ is the corresponding modular transformation.
%(see Lemma \ref{monodromy_transformation}).
\end{lemma}
\proof
By definition
\beq\label{sji}
(\Psi Re^{U/z}e_i,\d_j) = (-2\pi z)^{-3/2}\ z\d_j\, I_{\B_{i}}[e^{F/z}\omega],
\eeq
where $I_{\B_i}[e^{F/z}\omega]$ denotes the stationary phase asymptotics. Under analytic continuation the primitive form becomes $\omega/j(g,t_{-1})$, while the asymptotics $I_{\B_i}$ changes into $I_{P(\B_j)}$. It remains only to determine the monodromy of the flat vector fields
$$
\d_j=\sum_a (\d_js_a)\d/\d s_a\quad \mapsto \quad
\sum_{a} m_{aj}\,\d_a.
$$
Recall that under analytic continuation the flat coordinates
$t$ are transformed into $\nu(t)$ (see Lemma \ref{monodromy_transformation}). Put $s=(s_{-1},s_0,\dots,s_6)$. Then we have $s(\nu(t))=s(t)$. Using the chain rule
$$
\frac{Ds}{Dt}(\nu(t))\frac{D\nu}{Dt}(t) = \frac{Ds}{Dt}(t)
$$
we get that the matrix with entries $m_{aj}$ coincides with the Jacobian matrix $\Big(\frac{D\nu}{Dt}\Big)^{-1}$. The latter is straightforward to compute:
\beq\label{Jac}
\begin{bmatrix}
j(g,t_{-1})^2 & 0 & 0 & 0\\
              &   &   &  \\
-\frac{n_{12}^2}{2}\sum_{i=1}^6 t_it_{i'} & 1 & -\ge^{k} n_{12} m' & -\ge^{2k} n_{12}m''\\
              &   &   &  \\
j(g,t_{-1})n_{12}\leftexp{T}{m''} & 0 & \ge^k j(g,t_{-1})I_3 & 0 \\
              &   &   &  \\
j(g,t_{-1})n_{12}\leftexp{T}{m'} & 0 & 0 & \ge^{2k} j(g,t_{-1})I_3  \\
\end{bmatrix},
\eeq
where
$$
m'=[ t_6, t_5, t_4],\quad m''=[ t_3, t_2, t_1]
$$
and $\leftexp{T}{m'}$ and $\leftexp{T}{m''}$ are the columns with entries respectively (from top to bottom)
$ t_4, t_5, t_6$ and $ t_1, t_2, t_3$ Therefore, analytic continuation transforms the RHS of formula \eqref{sji} into
\beq\label{sji:1}
(-2\pi z)^{-3/2} \sum_{a=-1}^6 \, m_{aj}\, z\d_a\, \Big(I_{P(\B_i)}[ e^{F/z}\omega]/j(g,t_{-1})\Big).
\eeq
On the other hand the derivative in \eqref{sji:1} is
\ben
\delta_{a,-1} \frac{(-n_{12}z)}{j(g,t_{-1})^2} \, z\d_0\, I_{P(\B_i)}[e^{F/z}\omega] +
\frac{1}{j(g,t_{-1})}\,z\d_a\, I_{P(\B_i)}[e^{F/z}\omega].
\een
It is convenient to introduce the linear operator $P:\C^\mu\to \C^\mu$ whose action on the standard basis $\{e_i\}$ corresponds to the permutation of the Morse coordinate systems $\B_i\mapsto P(\B_i)$.  Then formula \eqref{sji:1} takes the form
\ben
\delta_{-1,j} (-n_{12}z) (\Psi R e^{U/z}P(e_i),\d_0)+
\sum_{a=-1}^6 j(g,t_{-1})^{-1}
( \Psi R e^{U/z}P(e_i) ,m_{aj} \d_a).
\een
Note that
\ben
\delta_{a,0}\delta_{-1,j} (-n_{12}z) + m_{aj}j(g,t_{-1})^{-1}
\een
is the $(a',j')$-entry of the matrix $M$  (see \eqref{matrix:M}). To finish the proof, it remains only to use that $\d_a=dt_{a'}$, $\d_j=dt_{j'}$.
\qed

\medskip

Define
\beq\label{J:matrix}
J(\nu,t_{-1})=
\begin{bmatrix}
1 & 0 \\
0 & j(g,t_{-1})^2
\end{bmatrix}
\oplus j(g,t_{-1})\, \ge^{2k}\,I_3 \oplus j(g,t_{-1})\, \ge^k\,I_3
\eeq
and
\beq\label{X:matrix}
X_{\nu,t_{-1}}(z)=
\begin{bmatrix}
1 & -n_{12}z/j(g,t_{-1}) \\
0 & 1
\end{bmatrix}\oplus I_6,
\eeq
for  $\nu=(g,k)\in {\rm SL}_2(\C)\times \Z.$
\begin{theorem}\label{desc:transf}
Analytic continuation transforms
\ben
\A_{t_{-1}}(\hbar;\q)\mapsto(\widehat{X}_{\nu,t_{-1}}{\A}_{t_{-1}})(\hbar j(\nu,t_{-1})^2;J(\nu,t_{-1})\q),
\een
where we first apply the operator $\widehat{X}_{\nu,t_{-1}}$ and then we rescale $\hbar$ and $\q$.
\end{theorem}
\proof
The idea is to  derive the transformation law for the ancestor potential $\A_t$ at some semi-simple point $t=(t_{-1},'t)$ and then pass to the limit $'t\to 0.$

According to Lemma \ref{ancestor:mon} the operator $\Psi R e^{U/z}$ is transformed into
$$\leftexp{T}{M}\Psi R e^{U/z}P.$$
We may assume that $P=1$ because $P$ is a permutation matrix, so its quantization $\widehat{P}$ will leave the product of Kontsevich--Witten tau functions invariant. Put $M=M_0+M_1 z.$ Then we have
\ben
\leftexp{T}{M}\Psi R e^{U/z} = \widetilde{\Psi} \widetilde{R} e^{U/z},\quad
\mbox{where}\quad \widetilde{\Psi}=M_0^{-1}\Psi,\quad \widetilde{R}=\Psi^{-1}M_0\leftexp{T}{M}\Psi R.
\een
The quantization is in general only a projective representation. However, the quantization of the operators $\Psi^{-1}M_0\leftexp{T}{M}\Psi$ and $R$ involves quantizing only $p^2$ and $p\,q$-terms. Since the cocycle \eqref{cocycle} on such terms vanishes we get
\ben
(\widetilde{R})\sphat = (\Psi^{-1}M_0\leftexp{T}{M}\Psi)\sphat\ \widehat{R}.
\een
The operators $M_0$ and $\Psi$ are independent of $z$ and their quantizations by definition are just changes of variables. Hence
\ben
(\widetilde{\Psi} \widetilde{R})\sphat = \widehat{M_0}^{-1} (M_0\leftexp{T}{M})\sphat \ (\Psi R)\sphat\ .
\een
By definition $\Delta_i^{-1}$ is $(\d_{u_i},\d_{u_i})$, which gains a factor of $j(\nu,t_{-1})^{-2}$ under analytic continuation. The ancestor potential \eqref{ancestor} is transformed into
\beq\label{ancestor:mon1}
\widehat{M_0}^{-1} (M_0\leftexp{T}{M})\sphat\ \Big( \A_t(j(\nu,t_{-1})^2\hbar ;j(\nu,t_{-1}) \q )\Big).
\eeq
Let us take the limit $'t\to 0$. Using formula \eqref{matrix:M}, we see that
\ben
M_0^{-1}\to j(\nu,t_{-1})\ J(\nu,t_{-1})^{-1},\quad
M_0\leftexp{T}{M}\to X_{\nu,t_{-1}}.
\een
It remains only to  notice that the rescaling
\ben
(\hbar,\q)\mapsto (j(\nu,t_{-1})^2\hbar,j(\nu,t_{-1})\q)
\een
commutes with the action of any quantized operator.
\qed

\subsection{Non-modular transformations}
Assume now that we have
two Frobenius structures corresponding to some cycles $A_1$ and $A_2$. According to Lemma \ref{flat:coord} the relation between the flat coordinates for $A_1$ and $A_2$ can be described by a
pair $\nu=(g,k)\in {\rm SL}_2(\C)\times \Z$.
\begin{lemma}\label{transf:general}
Let
$(\Psi R e^{U/z})_{A_i}$ be the asymptotic operator corresponding to the cycle $A_i(i=1,2);$ then
$(\Psi R e^{U/z})_{A_2}=\leftexp{T}{M}_\nu(\Psi R e^{U/z})_{A_1}.$
\end{lemma}
The proof of this Lemma is similar to the proof of Lemma \ref{ancestor:mon}. Moreover, using the same argument as in the proof of Theorem \ref{desc:transf}, we get
\begin{theorem}\label{transf:nonmod}
Let $\A_{A_i,t_{-1}}(\hbar;\q)$ be the total ancestor potential of the Frobenius structure corresponding to $A_i(i=1,2)$; then
\ben
\A_{A_2,\nu(t)_{-1}}(\hbar;\q)=(\widehat{X}_{\nu,t_{-1}}{\A}_{A_1,t_{-1}})(\hbar j(\nu,t_{-1})^2;J(\nu,t_{-1})\q),
\een
where we first apply the operator $\widehat{X}_{\nu,t_{-1}}$ and then we rescale $\hbar$ and $\q$.\end{theorem}
Let us emphasize the difference between Theorem \ref{desc:transf} and Theorem \ref{transf:nonmod}.
The former compares the values of the ancestor potential
at {\em  two different points} in $\mathbb{H}\times\C^{\mu-1}$. The latter compares the ancestor potentials of {\em two different} Frobenius structures at the {\em same point} in $\mathbb{H}\times\C^{\mu-1}$.

\subsection{The genus-1 potential}
We finish this section by describing the modular transformations of the genus-1 potential \eqref{genus1}.
The potential is a homogeneous function of degree 0 and therefore it depends only on the moduli $\tau=t_{-1}$ of the Frobenius structure. In fact, it is the derivative $\d F^{(1)}:=\d F^{(1)}/\d t_{-1}$ that transforms more naturally.
\begin{proposition}\label{transf:dgenus1}
Let $\nu=(g,k)\in {\rm SL}_2(\C)\times\Z$ be a modular transformation; then
\ben
\d F^{(1)}(\nu(t)) = j(g,t_{-1})^2 \d F^{(1)}(t)  +\Big(\frac{\mu}{24}-\frac{1}{2}\Big)\, n_{12}\, j(g,t_{-1}).
\een
\end{proposition}
\proof
Recall the notation in the proof of Theorem \ref{desc:transf}. Note that
\ben
\widetilde{R}_1=R_1+\Psi^{-1}M_0(\leftexp{T}{M}_1) \Psi
\een
and
\ben
(\d_{-1}\log \Delta_i)(\nu(t)) = j(g,t_{-1})^2(\d_{-1}\log \Delta_i)(t) + 2n_{12}j(g,t_{-1}).
\een
Since
\ben
\d F^{(1)}(\nu(t)) = \d_{-1}\Big( F^{(1)}(\nu(t))\Big) \ j(g,t_{-1})^2,
\een
in order to prove the proposition, we need to verify that
\ben
{\rm tr}\Big( \Psi^{-1}M_0(\leftexp{T}{M}_1) \Psi(\d_{-1}U)\Big) = -n_{12}/j(g,t_{-1}).
\een
On the other hand, since by definition
\ben
\Psi\, dU\, \Psi^{-1} = A,\quad A=\sum_{i=-1}^{\mu-2} (\d_i\bullet_t)\ dt_i,
\een
the LHS of the above identity is precisely
${\rm tr}\Big( M_0(\leftexp{T}{M}_1) (\d_{-1}\bullet_t)\Big).$ Since $F^{(1)}$ depends only on the moduli $t_{-1}$, we may assume that $'t=0.$ In this case however, the quantum multiplication operator $\d_{-1}\bullet_{t_{-1}}$ and $M_0(\leftexp{T}{M}_1)$ are given by the following matrices:
\ben
\begin{bmatrix}
0 & 0 \\
1& 0
\end{bmatrix}
\oplus (0\cdot I_6),\quad \mbox{and}\quad
\begin{bmatrix}
0 &  -n_{12}/j(g,t_{-1})\\
1& 0
\end{bmatrix}
\oplus (0\cdot I_6).
\een
The proposition follows.
\qed
\begin{corollary} \label{transf:genus1}
Let $\nu=(g,k)\in {\rm SL}_2(\C)\times\Z$ be a modular transformation; then
\ben
F^{(1)}(\nu(t)) = F^{(1)}(t) + \Big(\frac{\mu}{24}-\frac{1}{2}\Big)\, \log \, j(g,t_{-1}).
\een
\end{corollary}
Assume that we have two Frobenius structures corresponding to some cycles $A_i(i=1,2).$
 Let $F^{(1)}_{A_i}$ be the genus-1 potentials and $\nu=(g,k)\in {\rm SL}_2(\C)\times\Z$ be the transformation identifying the two flat structures; then
\begin{corollary}\label{transf:genus1b}
The following formula holds:
\ben
F_{A_2}^{(1)}(t) = F_{A_1}^{(1)}(t) + \Big(\frac{\mu}{24}-\frac{1}{2}\Big)\, \log \, j(g,t_{-1}).
\een
\end{corollary}
The genus-1 potential of a simple elliptic singularity was computed by I. Strachan \cite{Str}. In the $P_8$ case (when $\mu=8$) the answer is the following:
\beq\label{genus1:f}
F_A^{(1)}=-\frac{1}{24}\ \log \Big((27+\si^3)\pi_A^4\Big).
\eeq
The computation in \cite{Str} is carried out for a specific choice of the cycle $A$. However, using Corollary \ref{transf:genus1b}, we get that the above formula is valid for all possible choices of $A$.

\section{Anti-holomorphic completion}\label{sec:5}

  The transformation of the ancestor potential under the modular group from the previous section is
slightly complicated.
  In particular, the potential is not modular. A magic trick to restore the modularity is to complete it to
  an anti-holomorphic function. We call it an {\em anti-holomorphic completion}. It is motivated by physics
(see \cite{ABK}) and it has its origin
  in the so-called {\em holomorphic anomaly equations}. Mathematically, it can be thought as generalizing
  the construction of quasi-modular forms.

\subsection{Quasi-modular forms}
A function $f:\mathbb{H}\to \C$ is called a {\em holomorphic quasi-modular form of weight $m$} with
respect to some finite-index subgroup $\Gamma\subset {\rm SL}_2(\Z)$ if there
are functions $f_i$, $1\leq i\leq N$, holomorphic on $\mathbb{H},$ such
that
\begin{enumerate}
\item[(1)] The functions $f_0:=f$ and $f_i$ are holomorphic at the cusps of $\Gamma$.
\item[(2)] $$
f(\tau, \bar{\tau}) = f_0(\tau) +
f_1(\tau)(\tau-\overline{\tau})^{-1}+\dots
+f_N(\tau)(\tau-\overline{\tau})^{-N}.
$$
is modular, i.e.,
\ben
f(g\tau, g\overline{\tau}) = j(g,\tau)^m f(\tau, \overline{\tau}),\quad
\mbox{for all } g\in \Gamma,
\een
\end{enumerate}
$f(\tau, \overline{\tau})$ is called the {\em anti-holomorphic completion} of $f(\tau)$.

For more details we refer to \cite{KZ}. The key example of a
quasi-modular form is the Eisenstein series
\ben
G_2=-\frac{1}{24}+\sum_{n=1}^\infty \ \Big(\sum_{d|n} d\Big)\,
q^n,\quad q=e^{2\pi i\tau}.
\een
It is known that $G_2$ satisfies the following identity:
\ben
G_2(g(\tau)) = j(g,\tau)^2G_2(\tau)-\frac{1}{4\pi i}n_{12}
j(g,\tau),\quad g\in {\rm SL}_2(\Z),
\een
where the matrix $g$ and  its action on $\tau$ are the same as in
\eqref{action_on_tau}.   The map $\tau\mapsto g(\tau)$ induces the
following transformation:
\beq\label{quasi_modularity}
-\,\frac{1}{\tau-\overline{\tau}}\mapsto
-j(g,\tau)^2\frac{1}{\tau-\overline{\tau}} + n_{12}j(g,\tau).
\eeq
It follows that $G_2(\tau)-\frac{1}{4\pi i}(\tau-\overline{\tau})^{-1}$ transforms as
a modular form of weight 2.
It is not hard to prove that every quasi-modular form can be written
uniquely as a polynomial in $G_2$ whose coefficients are modular forms
on $\Gamma$.

For our purposes, we have to relax condition (1) in the definition of a quasi-modular forms. Namely, we will be assuming that the forms have finite order poles at the cusps.

\subsection{The anti-holomorphic completion of $\A_{t_{-1}}(\hbar;\q)$}

We continue to denote by $t\mapsto \nu(t)$ the transformation of the flat
coordinates corresponding to analytic continuation along a closed
loop in $\Si$. Recall also that $t_{-1}$ is identified with a point
$\tau'$ on the upper-half plane via some fractional linear
transformation $g$. Slightly abusing notation we
define complex conjugation:
\ben
t_{-1}=g(\tau'):=\frac{a\tau'+b}{c\tau'+d}\quad \mapsto \quad \overline{t}_{-1} =  \frac{a\overline{\tau'}+b}{c\overline{\tau'}+d},
\een
i.e., the complex conjugation of $t_{-1}$ is the one induced from the
upper half-plane. Note that the transformation law \eqref{quasi_modularity} still holds.
Put
\ben
\widetilde{X}_{t_{-1},\overline{t}_{-1}}(z) =
\begin{bmatrix}
1 & -z(t_{-1}-\overline{t}_{-1})^{-1} \\
&\\
0 & 1
\end{bmatrix}
\oplus I_6
\een
and define
\beq\label{desc:mod2}
\A_{t_{-1}, \overline{t}_{-1} }(\hbar;\q) =(\widetilde{X}_{t_{-1},\overline{t}_{-1}})\sphat \
\A_{t_{-1}}(\hbar;\q)
\ .
\eeq
As a consequence of Theorem \ref{desc:transf} we get the following corollary..
\begin{corollary}\label{c1}
Analytic continuation transforms the anti-holomorphic completion \eqref{desc:mod2} as follows:
\ben
\A_{t_{-1}, \overline{t}_{-1} }(\hbar;\q)\mapsto
\A_{t_{-1}, \overline{t}_{-1} }(j(\nu,t_{-1})^2\hbar;J(\nu,t)\q ).
\een
\end{corollary}
\proof
Using the transformation rule  \eqref{quasi_modularity}  we get that analytic continuation transforms
\ben
\widetilde{X}_{t_{-1},\overline{t}_{-1}}(z)\mapsto
\widetilde{X}_{t_{-1},\overline{t}_{-1}}(j(\nu,t_{-1})^2z)\, X_{\nu,t_{-1}}^{-1}(j(\nu,t_{-1})^2z).
\een
Since the quantization is a representation when restricted to the space of upper-triangular
symplectic transformations, the quantization of the RHS of the above equation is just a composition of the quantizations of the two operators. It remains only to use Theorem \ref{desc:transf} and the fact that under the rescaling
\ben
(\hbar,\q)\mapsto (j(\nu,t_{-1})^2\hbar,J(\nu,t_{-1})\q)
\een
the operators change as follows
\ben
(\widetilde{X}_{t_{-1},\overline{t}_{-1}}(z)\, X_{\nu,t_{-1}}^{-1}(z))\sphat\mapsto
(\widetilde{X}_{t_{-1},\overline{t}_{-1}}(j(\nu,t_{-1})^2z)X_{\nu,t_{-1}}^{-1}(j(\nu,t_{-1})^2z))\sphat.
\qed
\een

\subsection{Quasi-modularity of the ancestor potential}

It is convenient to introduce the following multi-index convention. Given $I=(i_0,i_1,\dots)$, with only finitely many $i_k:=(i_{k,-1},i_{k,0},\dots,i_{k,6})\in \Z^\mu$ different from $0$, we define the monomial
\ben
q_0^{i_0}(q_1+\one)^{i_1} q_2^{i_2}\cdots,
\een
where the raising of a vector variable to a vector power means raising
each component of the variable by the corresponding component of the
power and then taking their product. The anti-holomorphic ancestor potential
\eqref{desc:mod2} has the form
\ben
\A_{t_{-1}, \overline{t}_{-1}}(\hbar;\q)=\exp \sum_g \hbar^{g-1}
\F_{t_{-1}, \overline{t}_{-1}}^{(g)}(\q),
\een
where the genus-$g$ potential $\F_{t_{-1}, \overline{t}_{-1}}^{(g)}$ is a formal series of the following type:
\beq\label{genus:g}
\sum_I a^{(g)}_I(t_{-1}, \overline{t}_{-1})\ q_{0}^{i_{0}} (q_1+\one)^{i_1}q_2^{i_2}\cdots \ .
\eeq
Similarly, we let $a_I^{(g)}(t_{-1})$ be the coefficients of the ancestor potential $\A_{t_{-1}}(\hbar;\q).$
Finally, for each multi-index $I$ we introduce the following two integers:
\ben
d(I):=\sum_{k} (i_{k,-1}d_{-1}+\dots+i_{k,6}d_6),
\een
and
\ben
m(I):= \sum_{k} (2 i_{k,-1}+i_{k,1}+\dots+i_{k,6}).
\een
\begin{theorem}\label{quasi_mod}
The coefficient $a_I^{(g)}(t_{-1})$ is non-zero only if $d(I)$ is
an integer. Moreover,   each non-zero coefficient is a
quasi-modular form of weight $2g-2+m(I)$.
\end{theorem}
\proof
Analytic continuation transforms the series \eqref{genus:g}  into
\ben
\sum_I a^{(g)}_I(\nu(t_{-1}, \overline{t}_{-1}))\ q_{0}^{i_{0}}
(q_1+\one)^{i_1}q_2^{i_2}\cdots \ .
\een
On the other hand the substitution
\ben
\q\mapsto J(\nu,t)\q
\een
transforms the series \eqref{genus:g} into
\ben
\sum_I j(\nu,t_{-1})^{m(I)}e^{-2\pi id(I)} a^{(g)}_I(t_{-1}, \overline{t}_{-1})\
q_{0}^{i_{0}} (q_1+\one)^{i_1}q_2^{i_2}\cdots \ .
\een
Recalling Corollary \ref{c1}, we get the following formula:
\beq\label{modularity}
\widetilde{a}_I^{(g)}(\nu(t_{-1})) =
j(\nu,t_{-1})^{2g-2+m(I)}\ e^{-2\pi id(I)k}\ a_I^{(g)}(t_{-1}, \overline{t}_{-1}).
\eeq
Let us prove that $d(I)$ is an integer.  We claim that rescaling the asymptotics $\Psi R e^{U/z}$ via
\beq\label{rescaling}
t_i\mapsto e^{2\pi \sqrt{-1}\, d_i}t_i,\quad -1\leq i\leq 6,
\eeq
transforms  $\Psi R e^{U/z}$ according to the formula in
Lemma \ref{ancestor:mon}, where
$P=I_\mu$ and $\nu=(I_2,1).$
Assuming this claim, we get that under the rescaling
\eqref{rescaling} the descendant potential transforms according to
Theorem \ref{desc:transf}. This means that formula \eqref{modularity}
is still valid, and since $g=I_2$ we get that $e^{-2\pi i d(I)}=1$,
  i.e., $d(I)$ must be an integer.

The claim follows easily from the homogeneity property of the
oscillatory integrals. Namely, from
\ben
(z\d_z+E) \, J_{\B_i}(t,z)=
\theta\, J_{\B_i}(t,z),
\een
we get
\ben
J_{\B_i}(e^ct,e^cz)=e^{c\,\theta} \, J_{\B_i}(t,z),
\een
where the scalar $e^c$ acts on the $i$-th coordinate of $t$ with weight $d_i$. Let $c\to 2\pi i$
and note that the limit of the LHS, up to sign, coincides with rescaling $J_{\B_i}(t,z)$ via \eqref{rescaling}.
On the RHS the operator $e^{2\pi i\theta}$ coincides with the matrix $(-\leftexp{T}{M}_\nu),$ where $\nu=(I_2,1)$ (and $'t=0$).

\section{Relating to Gromov-Witten theory}\label{PXJ}
We have finished the proof of quasi-modularity of global Saito-Givental theory for $P_8$.
The remaining two cases of simple elliptic singularities are
\ben
X_9:\quad f(x,\si)=x_0^2x_2+x_0x_1^3+x_2^2+\si x_0x_1x_2\ ,\quad \si^3+27\neq 0,
\een
and
\ben
J_{10}:\quad f(x,\si)=x_0^3x_2+x_1^3+x_2^2+\si x_0x_1x_2\ ,\quad \si^3+27\neq 0.
\een
The proof of quasi-modularity for $X_9, J_{10}$ is
identical and we leave it for the readers to fill in the details. On the other hand, global Saito-Givental theory is considered to be a B-model theory.
To draw consequences for A-model theory, such as the GW theory of an elliptic orbifold $\P^1$, we have to identify the A-model
theory as a certain limit of the $B$-model. For our purposes, there are two important limits, the Gepner limit $\sigma=0$ and
the large complex structure limit $\sigma=\infty$. The Gepner limit corresponds to FJRW theory, while the large complex structure limit
corresponds to GW theory. For simple elliptic singularities, the appropriate flat coordinates at the Gepner limit $\sigma=0$
have been worked out already by Noumi-Yamada \cite{NY}. In our set-up, they correspond to the choice of a cycle $A$ such that
$\pi_A(0)=1, \pi'_A(0)=0$.
We define
\ben
\A_{Gepner, t}(\hbar;\q):=e^{F_A^{(1)}(t)+\frac{1}{24}\,t_{-1}}\ \A_{A,t}(\hbar;\q),
\een
where $F^{(1)}_A$ and $\A_{A,t}$ are the genus-1 potential and the total ancestor potential of the Frobenius structure corresponding to the cycles $A$. The following theorem is the LG-to-LG all genera mirror theorem
of Krawitz-Shen \cite{KS}.
\begin{theorem}
For $P_8, X_9, J_{10}$, $e^{-t_{-1}/24}\A_{Gepner, t}$ coincides with the ancestor potential function of
FJRW invariants, up to a linear identification of  the flat coordinates.
\end{theorem}
Since the FJRW ancestor potential function extends over the caustic, Lemma 4.3 implies
\begin{corollary}
The conjecture 3.2 holds for $P_8, X_9, J_{10}$, i.e., global Saito-Givental theory extends over the caustic.
\end{corollary}
The above corollary allows us to define the ancestor potential
function at $t_i=0$ for $i\geq 0$, which is crucial for our discussion of modularity.

\subsection{The divisor equation in singularity theory}
 Before we start to discuss the large complex structure limit and the relation to Gromov-Witten theory,
 we discuss the divisor equation in the B-model. The latter is an important tool in the computation of Gromov-Witten theory and it is necessary for the LG-to-CY mirror theorem of Krawitz--Shen \cite{KS}.

Let us denote by $P$ the flat vector field $\d/\d t_{-1}$. Let $t=(t_{-1},t_0,\dots,t_6)$ be a generic semi-simple
point. We will prove below that the correlators of the ancestor
potential $\A_t(\hbar;\q)$ are invariant with respect to
the transformation $t_{-1}\to t_{-1}+2\pi i$ and that they expand in a Fourier series in $q=e^{t_{-1}}$ near $q=0$ (see Proposition \ref{p8:prop2}).

\begin{comment}
Let us denote by $\A_{LCS,t}(\hbar;\q)$ the generating function whose correlators
are the Fourier expansions of the correlators of $\A_t(\hbar;\q)$. It satisfies the following differential equation:
\ben
q\d_q\ \A_{LCS,t}(\hbar;\q)=\d_{-1}\ \A_{LCS,t}(\hbar;\q),
\een
where $\d_{-1}:=\d/\d t_{-1}.$
\end{comment}

On the other hand, the ancestor potential of the singularity satisfies the differential equation
\ben
\d_{-1}\ \A_t = ((P\bullet_t/z)\sphat -\d_{-1} F^{(1)}(t))\A_t.
\een
This formula follows from the fact that the quantization operator $\Psi R e^{U/z}$ satisfies the quantum differential equations.

According to Corollary \ref{transf:dgenus1},  $\d_{-1}F^{(1)}(t)$ is a quasi-modular form (of weight 2). In particular, it  admits a Fourier expansion near $q=0$. Moreover, a straightforward computation shows that the constant term of this expansion is $-1/24.$ For $P_8$, one has to use formula \eqref{genus1:f} and the Fourier expansions of $\si$ and $\pi_A$ from Section \ref{lcs:p8}. In the other two cases the computation is again straightforward, thanks to the results of I. Strachan (see \cite{Str}). In other words,
\ben
\A_{LCS,t}(\hbar;\q):=e^{F^{(1)}(t)+\frac{1}{24}\,t_{-1} }\  \A_t(\hbar;\q)
\een
can be expanded into a Fourier series near $q=0$. Since the ancestor extends through the caustic, we can take the limit of $\A_{LCS,t}$  as $'t=(t_0,t_1,\dots,t_{\mu-2})\to 0.$ The resulting function, or more precisely its Fourier expansion, will be denoted by $\D_{LCS,q}(\hbar;\q).$ It satisfies the following differential equation:
\ben
q\d_q\ \D_{LCS,q}(\hbar;\q) =\Big((P\bullet /z)\sphat + \frac{1}{24}\Big) \D_{LCS,q}(\hbar;\q),
\een
where $P\bullet$ is the Frobenius multiplication by $P$ at $t=(t_{-1},0,\dots,0).$
Note that
$$
P\bullet\phi_i(x) = \delta_{i,0}P\quad \mbox{ for all } i.
$$
Therefore, the differential equation from above coincides with the divisor equation in the GW theory of the corresponding orbifold projective line with Novikov variable $q$ and divisor class $P.$
More precisely, the differential equation gives the following relation between the {\em correlators} of $\D_{LCS,q}(\hbar;\q)$:
\ben
\langle P,\phi_{a_1}\psi^k_1,\dots,\phi_{k_n}\psi^{k_n}\rangle_{g,n+1,d}
\een
equals
\ben
d \langle \phi_{a_1}\psi^k_1,\dots,\phi_{k_n}\psi^{k_n}\rangle_{g,n,d} +\sum_{i=1}^n \langle \dots,P\bullet \phi_i \psi_i^{k_i-1},\dots\rangle_{g,n,d},
\een
for all $(g,n,d)$, s.t., $d\neq 0$ or $2g-2+n>0$.
Here
\ben
\langle \phi_{a_1}\psi^k_1,\dots,\phi_{k_n}\psi^{k_n}\rangle_{g,n,d}
\een
is by definition the coefficient in front of $\hbar^{g-1}q_{k_1}^{a_1}\cdots q_{k_n}^{a_n}q^d$ in the generating function $\log\,(\D_{LCS,q})$.

\medskip

\begin{comment}
\subsection{The descendant potential in GW theory}
Let us point out that at the point $t=(\tau,0)$, the matrix $S_t$ is very simple, because the Frobenius multiplication satisfies

The differerential equation for $S_\tau$ with respect to $\tau$ is easily solved and it yields
\ben
S_\tau(z) =
 \begin{bmatrix}
1 & 0\\
\tau/z& 1
\end{bmatrix}
\oplus I_6.
\een
Therefore, it is straightforward to obtain the transformation law for the descendant potential as well.
\end{comment}

In the rest of this section, we focus on the coordinates at the large complex structure limit $\sigma=\infty$ and the identification
with GW theory. There is a large body of literature in the Calabi-Yau case. The identification is referred as {\em a mirror map}. In our case it goes
as follows. On the A-model side, the A-model moduli space has a coordinate $t_{-1}$ corresponding to a complexified K\"ahler class. GW theory involves power series in $q=e^{t_{-1}}$.
The B-model moduli are parameterized by $\sigma^3$ (for $P_8$ and $J_{10}$) or $\si^2$ for $X_9$. For our
purposes, it is convenient to work with $\sigma$ directly.
The mirror map is a map $\sigma\rightarrow t_{-1}(\sigma)$ defined locally around $\sigma=\infty$. By fixing a symplectic basis $\{A,B\}$ we can identify the K\"ahler class $t_{-1}$ with the modulus of the complex structure $\tau=\pi_B/\pi_A$.  We can treat the mirror map as a map $\tau \rightarrow  \tau(\sigma)$. The latter can be described explicitly via the Picard--Fuchs equations satisfied by the periods $\pi_A$ and $\pi_B$. We begin with the $P_8$-case.

\subsection{Large complex structure limit of the family $P_8$}\label{lcs:p8}
To begin with, let us construct a Frobenius manifold isomorphism between the Milnor ring
of $P_8$ and the quantum cohomology of $\mathbb{P}^1(3,3,3)$.

The orbifold cohomology admits the following natural basis:
$\Delta_0:=1$, $\Delta_{-1}:=P$ is the hyperplane class, and
 $\Delta_i$ and $\Delta_{i'}$, $i'=7-i$, are the cohomology classes 1
supported on the {\em twisted sectors} of the $i$-th orbifold point ($i=1,2,3$) of {\em
  age} $1/3$ and $2/3$ respectively (see \cite{CR} for some background on orbifold GW
theory). The only non-zero Poincar\'e pairings between these cohomology classes are
\ben
(\Delta_{-1},\Delta_0)=1,\quad (\Delta_i,\Delta_{i'})=1/3,\ i=1,2,\dots,6.
\een
According to Krawitz--Shen \cite{KS} the quantum cohomology and the higher-genus theory
are uniquely determined from the divisor equation and the correlators:
\ben
\langle \Delta_1,\Delta_2,\Delta_3\rangle_{0,3,1}=1,\quad \langle\Delta_i,\Delta_i,\Delta_i\rangle_{0,3,0}=1/3.
\een
To set up the B-model coordinates, we first have to choose a symplectic basis $\{A,B\}$.
The corresponding periods $\pi_A$ and $\pi_B$ are solutions to the differential equation
\eqref{Picard-Fuchs}, which has a regular singular point at $\si=\infty$. We choose the cycles
$A,B$ in such a way that
\beq\label{hgs}
\pi_{A}(\si)=-(-1)^{1/2}\,\si^{-1}\ \leftbase{2}{F}_1(1/3,2/3;1;-27/\si^3)..
\eeq
and
\beq\label{hgs:2}
\pi_{B}(\si)= -\frac{3}{2\pi i}\pi_{A}(\si) \log (-\si) +\frac{3}{2\pi i
  } \ (-3i\si^{-1})\sum_{k=1}^\infty b_k(-\si/3)^{-3k}\ .
\eeq
The coefficients $b_k$ can be determined uniquely from the recursion relation
\ben
-9k^2b_k+(9k^2-9k+2)b_{k-1} + (2k-1)a_{k-1}-2k a_{k}=0,
\een
where $a_k$ are the coefficients of the hypergeometric series $\leftbase{2}{F}_1(1/3,2/3;1;y)$, i.e.,
\ben
a_0=1,\quad a_k=\frac{(1/3)_k(2/3)_k}{(k!)^2},\ k\geq 1,\quad\mbox{where}\quad (b)_k=b(b+1)\cdots(b+k-1).
\een
These formulas suggest that the correct parameter of the B-model moduli is $(-\si/3)^{-3}$, but this is
not important for us. Exponentiating, we obtain an identity of the following form:
\ben
e^{2\pi i\tau/3} = -\si^{-1}\Big(1+\sum_{k=1}^\infty\  c_k\  (-\si)^{-3k}\Big),
\een
where the coefficients $c_k$ can be written explicitly in terms of $a_k$ and $b_k$. By inverting
the above series we can obtain the Fourier expansion of $-\si^{-1}$ in terms of $q:=e^{2\pi i\tau/3}.$
Note that
\ben
{\rm res}_{x=0}\ \frac{x_0x_1x_2}{f_{x_0}f_{x_1}f_{x_2}}\,d^3x = \frac{1}{\si^3+27}
\een
and that the monomials $(\si^3+27)^{\rm deg\, \phi_i}\phi_i(x)(-1\leq i\leq 6)$ provide a basis in which the residue pairing with respect to the standard volume form is constant. It follows that the identifications
\ben
1& =& 1\\
P& =&  (\si^3+27)\pi_A^2(\si)\phi_{-1}(x)\\
27^{\rm deg \, \Delta_i-\frac{1}{3}}\Delta_i &= & (-1)^{{\rm deg\, \phi_i}-\frac{1}{2} } (\si^3+27)^{\rm deg\, \phi_i}\phi_i(x)\, \pi_A(\si) ,\quad
1\leq i \leq 6,
\een
provide an isomorphism between the Poincar\'e and the residue pairings. Moreover, after a straightforward
computation, we find the leading terms of the Fourier series of the following correlators:
\ben
\langle \Delta_1,\Delta_2,\Delta_3\rangle_{0,3}=-i\pi_A(\si)=q+q^4+2q^7+2q^{13}+q^{16}+2q^{19}+\cdots,
\een
\ben
\langle \Delta_1,\Delta_1,\Delta_1\rangle_{0,3} =-i (-\si/3)\pi_A(\si)=\frac{1}{3}+2q^3+2q^9+2q^{12}+\cdots\ .
\een
\begin{remark}
The Fourier expansion of $-i\pi_A(\si)$ coincides with Saito's eta product (see \cite{Sa4})
\ben
\eta_{E_6^{(1,1)}}(3\tau):=\eta(9\tau)^3\eta(3\tau)^{-1}.
\een
\end{remark}
\begin{proposition}\label{p8:prop1}
The cycles $A$ and $B$ are integral up to a scalar factor and $\tau$ is a modulus of the elliptic curve at infinity.
\end{proposition}
\proof
The $j$-invariant  of the elliptic curve at infinity is
\ben
j(\si)= -\frac{\si^3(-216+\si^3)^3}{(27+\si^3)^3}.
\een
According to Kodaira \cite{Ko}, there exists a symplectic basis $\{A',B'\}$ of $H_1(E_{\si};\Z)$ whose
monodromy around $\si=\infty$ is the same as the monodromy of $\{\pi_A,\pi_B\}$. This implies that
\ben
A=c\, A',\quad B=c\,B'+ d\,A'
\een
for some constants $c$ and $d$. On the other hand the Fourier expansion of the above $j$-invariant is
\ben
\frac{1}{q^3}+744+196884 q^3+21493760 q^6+\cdots,\quad q=e^{2\pi i\tau/3}.
\een
Comapring with the well known Fourier expansion of the $j$-invariant, we get that $\tau=\tau'$
and hence the constant $d=0$. \qed

Note that under the identification between the quantum cohomology and the Milnor ring from above,
the K\"ahler parameter (i.e., the coordinate along the hyperplane class $P$) becomes $t_{-1}$.  The next
Proposition guarantees that the divisor equations in singularity theory and in Gromov--Witten theory are
the same.
\begin{proposition}\label{p8:prop2}
The K\"ahler parameter is related to $\tau$ via the following mirror map: $t_{-1}=2\pi i\tau/3.$
\end{proposition}
\proof
We want to compute the constant $(1,\d/\d\tau)_A$ (the index $A$ means residue pairing with respect to
$d^3x/\pi_A$). To begin with note that
\beq\label{div:1}
\frac{\d \tau}{\d \si}\ (1,\d/\d\tau)_A = (1,\d/\d \si)_A = (1,x_0x_1x_2)_A = \frac{1}{(27+\si^3)\pi_A^2}.
\eeq
On the other hand
\ben
\frac{\d \tau}{\d \si} = \frac{\pi_B'\pi_A-\pi_B\pi_A'}{\pi_A^2}.
\een
The numerator is the Wronskian of the solutions $\pi_B$ and $\pi_A$ of the differential equation
\eqref{Picard-Fuchs} and hence it equals
\ben
{\rm const}\cdot (27+\si^3)^{-1} \quad \sim\quad {\rm const}\cdot \si^{-3},
\een
where we took the expansion near $\si=\infty$ and kept only the leading term. On the other hand
using the expansions
of $\pi_A$ and $\pi_B$ at $\si=\infty$ (see formulas \eqref{hgs} and \eqref{hgs:2}) one can check that the leading order term of the numerator is: $3\si^{-3}/2\pi i$. Therefore, the above constant is $3/2\pi i.$
Now, from equation \eqref{div:1}, we get $(1,\d/\d \tau)_A=2\pi i/3.$ On the other hand, since
$1=(1,P)=(1,\d/\d t_{-1})$, we must have $t_{-1}=2\pi i\tau/3.$
\qed

All necessary conditions for the reconstruction theorem of Krawitz--Shen are satisfied. Therefore, we have the following theorem:
 \begin{theorem}
 Under the identification of $q=e^{2\pi i \tau/3}$ with the Novikov variable, $\D_{LCS,q}$ is equal to the descendant potential function of the elliptic orbifold $\P^1$ with weights $(3,3,3)$.
\end{theorem}
 It is well known that the modular group, i.e., the monodromy group of the local system
$H_1(E_\si;\Z),\si\in \Si,$ is $\Gamma(3)$ -- the principal congruence subgroup of level 3. Let us denote
by $a_I^{(g)}(\tau)$ the coefficients  of the ancestor potential $\A_{t_{-1}}$ of the singularity. An immediate consequence is
 \begin{corollary}\label{cor:p8}
The following statements hold:
 \begin{itemize}
 \item[(1)] $a^{(g)}_I(\tau)$ has no pole at the cusp $\tau=i\,\infty$.
 \item[(2)] $\A_{Gepner}$ is related to $\A_{LCS}$ by the composition of analytic continuation and the quantization of a symplectic transformation (see Theorem \ref{transf:nonmod}).
 \item[(2)] The coefficients of the total descendant potential for the GW theory of $\P^1(3,3,3)$ are quasi-modular forms on $\Gamma(3)$.
 \end{itemize}
 \end{corollary}
 \begin{remark}
 From the B-model alone, it is difficult to see whether $a_I^{(g)}(\tau)$ does not have a pole at the cusp, i.e., at $q=0$. The situation is similar to the extendibility of Givental's function to the caustic. Again we draw the conclusion from the mirror A-model side by using Krawitz-Shen's GW-to-LG all genera mirror theorem.
\end{remark}
\begin{remark}
One may wonder if analytic continuation alone will relate $\A_{Gepner}$ to $\A_{LCS}$. The answer is generally no.
We are choosing different symplectic bases at $\sigma=0, \sigma=\infty$. One basis may not be analytic
continuation of the other. For example, we can often choose an integral basis at $\sigma=\infty$. But the basis at $\sigma=0$
is not integral in general.
\end{remark}
\begin{remark}
It is convenient to use the language of symplectic bases to describe the ideas. Technically, it is easier to work with Picard-Fuchs equations.
Fortunately, the two approaches are equivalent. However, it is  generally a difficult question to identify a
solution of the Picard-Fuchs equation with the period of an explicit cycle.
\end{remark}

\subsection{Large complex structure limit of the family $X_9$}
The primitve forms are given by $d^3x/\pi_A(\si),$ where $\pi_A(\si)$ is a solution to the same differenrial equation as in the $P_8$-case.
\begin{comment}
of the differential equation
\beq\label{PF:X9}
\frac{d^2u}{d\si^2}+\frac{2\si}{\si^2-4}\,
\frac{du}{d\si}+\frac{1}{4(\si^2-4)}\, u=0.
\eeq
\end{comment}
\subsubsection{The Gauss--Manin connection in the marginal direction}
In order to identify the quantum cohomology with the Milnor ring, we need to
find the monomials in the Milnor ring for which the residue pairing assumes a constant form. We fix the following basis in the Milnor ring: $\phi_{-1}=x_0x_1x_2$, $\phi_0=1$, and $\phi_i$ for $i=1,2,\dots, 7$ are given respectively by
\ben
x_0,\quad x_1,\quad x_2,\quad x_1^2,\quad
x_0x_1,\quad x_0x_2,\quad x_1x_2.
\een
Put $\Phi_i(\si)$ for the section $\int\phi_i d^3x/df$ of the middle cohomology bundle. A straightforward computation, similar to the one in Section \ref{modular_transformations}, gives the following differential equations:
\ben
\d_\si \Phi_1 & = &-\frac{\si^2}{4(27+\si^3)}\, \Phi_1-\frac{9}{2(27+\si^3)}\, \Phi_2,\\
\d_\si\Phi_2 & = & \frac{3\si}{4(27+\si^3)}\, \Phi_1-\frac{\si^2}{2(27+\si^3)}\, \Phi_2,
\een
\ben
\d_\si \Phi_3 & = & -\frac{\si^2}{27+\si^3}\,\Phi_3-\frac{9}{2(27+\si^3)}\, \Phi_5\\
\d_\si \Phi_4 & = & -\frac{9}{27+\si^3}\,\Phi_3 +\frac{3\si}{2(27+\si^3)}\, \Phi_5,\\
\d_\si \Phi_5 & = & \frac{3\si}{27+\si^3}\, \Phi_3 -\frac{\si^2}{2(27+\si^3)}\, \Phi_5,
\een
and
\ben
\d_\si\Phi_6 & = & -\frac{7\si^2}{4(27+\si^3)}\, \Phi_6-\frac{9}{2(27+\si^3)}\, \Phi_7\\
\d_\si\Phi_7 & = & \frac{21\si}{4(27+\si^3)}\, \Phi_6-\frac{\si^2}{2(27+\si^3)}\, \Phi_7.
\een
From here we get the following solutions:
\ben
\Phi_1(\si)&=&\si^{-1/4}\Phi_{1,1}(\si)A_1+\si^{-5/2}\Phi_{1,2}(\si)A_2,\\
\Phi_2(\si)&=&-\si^{-5/4}\Phi_{2,1}(\si)A_1+\frac{1}{2}\si^{-1/2}\Phi_{2,2}(\si)A_2,
\een
\ben
\Phi_3(\si)&=&\si^{-1}\Phi_{3,1}(\si)A_3+\si^{-5/2}\Phi_{3,2}(\si)A_5,\\
\Phi_4(\si)&=&A_4-\frac{1}{3}\si\Phi_3(\si),\\
\Phi_5(\si)&=&-2\si^{-2}\Phi_{5,1}(\si)A_3+
\frac{1}{3}\si^{-1/2}\Phi_{5,2}(\si),
\een
\ben
\Phi_6(\si)&=&\si^{-7/4}\Phi_{6,1}(\si)A_6+\si^{-5/2}\Phi_{6,2}(\si)A_7\\
\Phi_7(\si)&=&-\frac{7}{3}\si^{-11/4}\Phi_{7,1}(\si)A_6+\frac{1}{6}\si^{-1/2}\Phi_{7,2}(\si)A_7,
\een
where
\ben
\Phi_{1,1}(\si)&=&
\leftbase{2}{F}_1(1/12,5/12;1/4;-27/\si^3),\\
\Phi_{1,2}(\si)&=&
\leftbase{2}{F}_1(5/6,7/6;7/4;-27/\si^3),
\een
\ben
\Phi_{2,1}(\si)&=&
\leftbase{2}{F}_1(5/12,13/12;5/4;-27/\si^3),\\
\Phi_{2,2}(\si)&=&
\leftbase{2}{F}_1(1/6,5/6;3/4;-27/\si^3),
\een
\ben
\Phi_{3,1}(\si)&=&
\leftbase{2}{F}_1(1/3,2/3;1/2;-27/\si^3),\\
\Phi_{3,2}(\si)&=&
\leftbase{2}{F}_1(5/6,7/6;3/2;-27/\si^3),
\een
\ben
\Phi_{5,1}(\si)&=&
\leftbase{2}{F}_1(2/3,4/3;3/2;-27/\si^3),\\
\Phi_{5,2}(\si)&=&
\leftbase{2}{F}_1(1/6,5/6;1/2;-27/\si^3),
\een
\ben
\Phi_{6,1}(\si)& =&
\leftbase{2}{F}_1(7/12,11/12;3/4;-27/\si^3),\\
\Phi_{6,2}(\si)&=&
\leftbase{2}{F}_1(5/6,7/6;5/4;-27/\si^3),
\een
\ben
\Phi_{7,1}(\si)&=&
\leftbase{2}{F}_1(11/12,19/12;7/4;-27/\si^3),\\
\Phi_{7,2}(\si)&=&
\leftbase{2}{F}_1(1/6,5/6;1/4;-27/\si^3),
\een
and $A_i$ are flat sections of the middle cohomology bundle. Solving for $A_i$ in terms of $\Phi_i$ we get certain polynomials in the Milnor ring which, according to our general construction of flat coordinates, should be part of a basis in which the residue pairing (with respect to the standard form $d^3x$) is constant.
\begin{comment}
\ben
(\si^2-4)^{q_i}\phi_i\pi_A,\quad i=1,2,4,6,7,
\een
and
\ben
(\si+2)^{1/2}(\phi_3+\phi_5)\pi_A\quad\mbox{and}\quad (\si-2)^{1/2}(\phi_3-\phi_5)\pi_A.
\een
\end{comment}
\subsubsection{The orbifold quantum cohomology}
The orbifold cohomology of $\P^1(4,4,2)$ has the following natural basis: 1 is the unit, $P$ is the hyperplane class,
and the remaining cohomology classes are supported on the twisted
sectors. Namely,
$\Delta_{i1},\Delta_{i2}, \Delta_{i3}$, $i=1,2$ are the units
($\Delta_{ik}$ has degree $k/4$) in the
cohomology of the twisted sectors of the $i$-th $\Z/4\Z$-orbifold
point, and  $\Delta_{31}$ is the unit in the cohomology of the twisted
sector of the $\Z/2\Z$-orbifold point. The only non-zero pairings are
\ben
(1,P)=1,\quad (\Delta_{i1},\Delta_{i3})=1/4,\quad
(\Delta_{i2},\Delta_{i2})=1/4,\quad (\Delta_{31},\Delta_{31})=1/2,
\een
where $i=1,2.$ Also, the following correlators are easily computed because the
corresponding moduli spaces are points.
\ben
\langle\Delta_{i1},\Delta_{i1},\Delta_{i2}\rangle_{0,3,0}=1/4,\quad i=1,2,
\een
\ben
\langle\Delta_{11},\Delta_{21},\Delta_{31}\rangle_{0,3,1}=1.
\een
According to Krawitz--Shen \cite{KS}, the quantum cohomology and the higher genus theory are uniquely determined from the above relations and the divisor equation.

We specify the cycles $\{A,B\}$ by choosing the corresponding periods $\pi_A$ and $\pi_B$. The period $\pi_A$
is the same as in the $P_8$-case (see \eqref{hgs}), while $\pi_B$ is 3 times larger:
\ben
\pi_{B}(\si)= -\frac{9}{2\pi i}\pi_{A}(\si) \log (-\si) +\frac{9}{2\pi
  i} \ (-3i\si^{-1})\sum_{k=1}^\infty b_k(-\si/3)^{-3k}\ .
\een
\begin{comment}
Let us construct the identification with the Milnor ring.
At $\si=\infty$ the equation has the following anti-invariant solution:
\ben
\pi_A(\si)=\si^{-1/2}\ \leftbase{2}{F}_1(1/4,3/4;1;4/\si^2).
\een
The differential equation admits a second solution of the following form:
\beq\label{2-nd_sol}
\pi_B(\si)=-\frac{2}{2\pi i}\, \pi_A(\si)\,
\log(-4\si)+\frac{1}{2\pi i}\,\si^{-1/2}\,\sum_{k=1}^\infty b_k\,(-2/\si)^{2k},
\eeq
where the coefficients $b_k$ are determined recursively from the relations
\ben
k^2 b_k= (k^2-k+3/16)b_{k-1} - 2ka_k+(2k-1)a_{k-1},\quad b_0=0,
\een
and $a_k$ are the coefficients of the hypergeometric series, i.e.,
\ben
a_0=1,\quad a_k=\frac{(1/4)_k (3/4)_k}{(k!)^2},\ k\geq 1, \mbox{  where } (b)_k = b(b+1)\cdots (b+k-1).
\een
First we divide \eqref{2-nd_sol} by $\pi_A$ and set $\tau:=\pi_B/\pi_A$. Then we multiply by $2\pi i/4$ and
exponentiate. We get
\ben
2\sqrt{2}\, q=(-2/\si)^{1/2}\,\Big( 1 + \sum_{k=1}^\infty c_k (-2/\si)^{2k}\Big),\quad q:=e^{2\pi i\tau/4}.
\een
From here one can easily express $(2/\si)^{1/2}$ as a power series in $q$ of the following form:
\ben
(-2/\si)^{1/2}=2\sqrt{2} \, q\Big(1+ \sum_{k=1}^\infty d_k q^{4k}\Big).
\een
\end{comment}
\begin{proposition}\label{x9:prop1}
Up to a scalar, the cycles $A$ and $B$ are integral and $\tau$ is a modulus of the elliptic curve at infinity.
\end{proposition}
\proof
The argument is similar to the one in Proposition \ref{p8:prop1}: we need to check that the $j$-invariant
has the correct Fourier expansion in terms of $e^{2\pi i\tau}.$
The $j$-invariant of the elliptic curve at infinity is
\ben
j(\si) = -\frac{(24\si+\si^4)^3}{(27+\si^3)}.
\een
Substituting in this formula the Fourier series of $-\si^{-1}$, we get
\ben
e^{-2\pi i\tau}+744+196884 e^{2\pi i\tau}+\cdots.\qed
\een
\begin{comment}
Put
\ben
\gamma=\frac{1}{2}\Big((\si-2)^{1/2}+(\si+2)^{1/2}\Big).
\een
\end{comment}
Since
\ben
{\rm res}_{x=0}\ \frac{x_0x_1x_2}{f_{x_0}f_{x_1}f_{x_2}}\, d^3x =
\frac{9}{4(27+\si^3)},
\een
it is easy to check that if we identify
\ben
1=1,\quad P=\frac{4}{9}(27+\si^3)\, x_0x_1x_2\,\pi_A^2,
\een
and
\ben
\Delta_{11}&=&
2e^{\pi\sqrt{-1}/4}(27+\si^3)^{1/4}\Big(
\frac{1}{2}\si^{-1/2}\Phi_{2,2}(\si)\, x_0-
\si^{-5/2}\Phi_{1,2}(\si)\, x_1\Big) \,\pi_A\\
\Delta_{12}&=&
-(27+\si^3)^{1/2}\Big(
2\si^{-2}\,\Phi_{5,1}(\si)\, x_2+
\si^{-1}\,\Phi_{3,1}(\si)\, x_1x_2\Big)\pi_A, \\
\Delta_{13}&=&
-2e^{-\pi\sqrt{-1}/4}(27+\si^3)^{3/4}\Big(
\frac{1}{6}\si^{-1/2}\Phi_{7,2}(\si)\, x_0x_2-
\si^{-5/2}\Phi_{6,2}(\si)\, x_1x_2\Big)\pi_A \\
\Delta_{21}&=&
(27+\si^3)^{1/4}\Big(
\si^{-1/4}\Phi_{1,1}(\si)\,x_1+\si^{-5/4}\Phi_{2,1}(\si)\, x_0\Big)\pi_A\\
\Delta_{22}&=& e^{\pi \sqrt{-1}/2}\Big(x_1^2+\frac{1}{3}\si x_2+ (27+\si^3)^{1/2}\Big(\frac{1}{3}\si^{-1/2}\Phi_{5,2}(\si)\,x_2-\si^{-5/2}\Phi_{3,2}(\si)\, x_0x_1\Big)\Big)
\pi_A \\
\Delta_{23}&=& (27+\si^3)^{3/4}\Big(\si^{-7/4}\Phi_{6,1}(\si)\,x_1x_2+\frac{7}{3}\si^{-11/4}\Phi_{7,1}\,x_0x_2\Big) \pi_A\\
\Delta_{31}&=&e^{\pi \sqrt{-1}/2}\Big(-x_1^2-\frac{1}{3}\si x_2+2(27+\si^3)^{1/2}\Big(\frac{1}{3}\si^{-1/2}\Phi_{5,2}(\si)\,x_2-\si^{-5/2}\Phi_{3,2}(\si)\, x_0x_1\Big)\Big)
\pi_A
\een
then the residue and the Poincar\'e pairings coincide.  Moreover, put $q=e^{2\pi i\tau/4}$; then we have the following formulas for the correlators:
\ben
\langle \Delta_{11},\Delta_{21},\Delta_{31}\rangle_{0,3}=
 q + 2 q^5 + q^9 + 2 q^{13} + 2 q^{17} + 3 q^{25} + O(q^{26}),
\een
\ben
\langle \Delta_{11},\Delta_{11},\Delta_{12}\rangle_{0,3}=
\frac{1}{4} + q^4 + q^8 + q^{16} + 2 q^{20} + O(q^{26}),
\een
 and
\ben
\langle \Delta_{11},\Delta_{11},\Delta_{22}\rangle_{0,3}=
q^2 + 2 q^{10} + q^{18} +O(q^{26}).
\een
\begin{proposition}\label{x9:prop2}
The K\"ahler parameter is related to $\tau$ via the following mirror map: $t_{-1}=2\pi i\tau/4.$
\end{proposition}
The proof is along the same lines as Proposition \ref{p8:prop2} and it is left as an exercise.
Using again the results of Krawitz-Shen we have the following GW-to-LG all genera mirror theorem.
\begin{theorem}\label{lcsl:x9}
Let $q=e^{2\pi i \tau/4}$. Under the above isomorphism between the quantum cohomology and the Milnor ring, $\D_{LCS,q}$ is equal to the total descendant potential function of the elliptic orbifold $\P^1$ with weights $(4,4,2)$.
\end{theorem}
\begin{comment}
Let $\Gamma(X_9)$ be the monodromy group of the elliptic fibration
\ben
\{X_0^3X_2+X_1^3+X_2^2+\si X_0X_1X_2=0\}\subset \mathbb{P}^2(1,1,2),\quad \si^2\neq 4.
\een
$\Gamma(X_9)$ is isomorphic to the subgroup of ${\rm SL}_2(\Z)$ consisting
of those matrices whose diagonal entries are $1$ modulo $4$ (see Appendix \ref{app:B}).
In other words, the total ancestor potential $A_{t_{-1}}$ for the $X_9$ singularity is quasi-modular on $\Gamma(X_9)$. On the other hand, in order to identify the correlators in GW theory with functions on the modular curve, we need to extend the space of rational functions by allowing also $(-\si+2)^{1/2}$ and $(-\si-2)^{1/2}$. Therefore, in GW theory we have modularity only with respect to the subgroup of $\Gamma(X_9)$ that leaves these two functions invariant. It is easy to see that this subgroup is the principal congruence subgroup $\Gamma(4).$
\end{comment}

Theorem \ref{quasi_mod} and \ref{lcsl:x9} yield

\begin{corollary}
The Gromov-Witten total descendant potential function of the elliptic orbifold $\P^1$ with weights $(4,4,2)$ is quasi-modular for $q=e^{2\pi i \tau/4}$ and for some finite index subgroup of ${\rm SL}_2(\Z)$.
\end{corollary}

\begin{comment}
In particular, using that the ring of modular forms for $\Gamma(4)$ is generated by $\theta_{00}^2,\theta_{01}^2,\theta_{10}^2$, where
\ben
\theta_{00}=\sum_{k\in \Z} q^{2k^2},
\quad
\theta_{01}=\sum_{k\in \Z} (-1)^kq^{2k^2},
\quad
\theta_{10}=2q^{1/2}\, \sum_{k=1}^\infty q^{2(k^2+k)},
\een
we get the following formulas:
\ben
2i\pi_A=\theta_{10}^2,\quad 2i\gamma\pi_A=\theta_{00}^2+\theta_{01}^2,\quad
2i\gamma^{-1}\pi_A=\theta_{00}^2-\theta_{01}^2.
\een
\end{comment}

\subsection{Large complex structure limit of  the family $J_{10}$}
The primitive forms are given by $d^3x/\pi_A(\si)$, where $\pi_A$ is a solution to the same differential equation as in the $P_8$-case.

\subsubsection{The Gauss--Manin connection in the marginal direction}
We fix the following basis in the Milnor ring: $\phi_{-1}=x_0x_1x_2$, $\phi_0=1$, and $\phi_i$ for $i=1,2,\dots,8$ are given by the monomials
\ben
x_0,\quad x_0^2,\quad x_1,\quad x_0^3,\quad x_0x_1,\quad x_0^4,\quad x_1^2,\quad x_0^5.
\een
Put $\Phi_i(\si)$ for the section $\int \phi_i(x)d^3x/df$ of the middle cohomology bundle. They satisfy the
following system of differential equations:
\ben
\d_\si\Phi_1 = -\frac{\si^2}{2(27+\si^3)}\, \Phi_1,\quad
\d_\si\Phi_8 = -\frac{3}{2(27+\si^3)}\Big(\frac{24(18+\si^3)}{24\si+\si^4}+\si^2\Big)\, \Phi_8,
\een
\ben
\d_\si\Phi_2=-\frac{9}{27+\si^3}\, \Phi_3,\quad \d_\si\Phi_3 = -\frac{\si^2}{27+\si^3}\, \Phi_3,
\een
\ben
\d_\si\Phi_4=\frac{\si^2}{2(27+\si^3)}\,\Phi_4-\frac{18}{27+\si^3}\, \Phi_5,\quad
\d_\si\Phi_5=-\frac{2\si^2}{27+\si^3}\,\Phi_5-\frac{3\si}{2(27+\si^3)}\, \Phi_4,
\een
\ben
\d_\si\Phi_6=\frac{1}{\si}\, \Phi_6-\frac{162}{\si^2(27+\si^3)}\,\Phi_7,\quad
\d_\si\Phi_7=\frac{27-\si^3}{\si(27+\si^3)}\,\Phi_7.
\een
The solutions have the form
\ben
\Phi_1(\si)=(27+\si^3)^{-1/6}A_1,\quad \Phi_8(\si)=\frac{24\si+\si^4}{(27+\si^3)^{5/6}}A_8
\een
\ben
\Phi_2(\si)=(-\si/3)(27+\si^3)^{-1/3}A_3+A_2,\quad \Phi_3(\si)=(27+\si^3)^{-1/3}A_3,
\een
\ben
\Phi_4(\si)=\Phi_{41}(\si)\,A_4+ \Phi_{42}(\si)\,A_5,\quad
\Phi_5(\si)=\Phi_{51}(\si)\,A_4+ \Phi_{52}(\si)\,A_5,
\een
\ben
\Phi_6(\si)=\frac{18+\si^3}{3(27+\si^3)^{2/3}}A_7 +\si A_6,\quad
\Phi_7(\si)=\si(27+\si^3)^{-2/3}A_7,
\een
where $A_i$ are flat sections of the middle cohomology bundle and
\ben
\Phi_{41}=&3(-\si/3)^{1/2}& \leftbase{2}{F}_1(-1/6,1/6;-1/2;-27/\si^3),\\
\Phi_{42}=&\frac{4}{9}(-\si/3)^{-4}& \leftbase{2}{F}_1(\ \ 4/3, 5/3;\ \ 5/2;-27/\si^3),\\
\Phi_{51}=&(-\si/3)^{-1/2}& \leftbase{2}{F}_1(\ \ 1/6,5/6;\ \ 1/2;-27/\si^3),\\
\Phi_{52}=&(-\si/3)^{-2}& \leftbase{2}{F}_1(\ \ 2/3,4/3;\ \ 3/2;-27/\si^3)..
\een
From here we can determine the elements in the Milnor ring that correspond
to the flat sections $A_i$. They should correspond to orbifold cohomology
classes of $\P^1(6,3,2)$.

\subsubsection{The quantum cohomology}
A natural basis in the orbifold cohomology is:
the unit 1, the hyperplane class $P$, and the units of the twisted sectors
$\Delta_{1i} (1\leq i\leq 5)$, $\Delta_{2j}(j=1,2)$, and $\Delta_{31}$.
The Poincar\'e pairing in this basis is non-zero only in the following cases:
\ben
(\Delta_{1,i},\Delta_{1,j})=\delta_{i+j,6}/6,\quad
(\Delta_{2,1},\Delta_{2,2})=1/3,\quad
(\Delta_{31},\Delta_{3,1})=1/2.
\een
According to Krawitz--Shen \cite{KS}, the quantum cohomology
and the higher-genus theory are uniquely determined by the divisor equation and the
following correlators:
\ben
\langle\Delta_{11},\Delta_{11},\Delta_{14}\rangle_{0,3,0}=1/6,&
\langle\Delta_{11},\Delta_{12},\Delta_{13}\rangle_{0,3,0}=1/6 \\
\langle\Delta_{21},\Delta_{21},\Delta_{21}\rangle_{0,3,0}=1/3, &
\langle\Delta_{11},\Delta_{21},\Delta_{31}\rangle_{0,3,1}=1.
\een
We specify the cycles $\{A,B\}$ by choosing the corresponding periods $\pi_A$ and $\pi_B$. The period $\pi_A$
is the same as in the $P_8$-case (see \eqref{hgs}), while $\pi_B$ is 3 times larger:
\ben
\pi_{B}(\si)= -\frac{9}{2\pi i}\pi_{A}(\si) \log (-\si) +\frac{9}{2\pi
  i} \ (-3i\si^{-1})\sum_{k=1}^\infty b_k(-\si/3)^{-3k}\ .
\een
From here we can express $-\si^{-1}$ as a Fourier series in $e^{2\pi i\tau/9}$, where $\tau:=\pi_B/\pi_A$.
\begin{proposition}\label{j10:prop1}
Up to a scalar, the cycles $A$ and $B$ are integral and $\tau$ is a modulus of
the elliptic curve at infinity.
\end{proposition}
\proof
The argument is similar to the one in Proposition \ref{p8:prop1}: we need to check that the $j$-invariant
has the correct Fourier expansion in terms of $e^{2\pi i\tau}.$
The $j$-invariant of the elliptic curve at infinity is
\ben
j(\si) = -\frac{(24\si+\si^4)^3}{(27+\si^3)}.
\een
Substituting in this formula the Fourier series of $-\si^{-1}$, we get
\ben
e^{-2\pi i\tau}+744+196884 e^{2\pi i\tau}+\cdots.\qed
\een
In order to match the quantum cohomology and the Milnor ring we make the following identifications
($\eta=e^{2\pi i/6}$):
\ben
&&(-1)^{\frac{1}{2}-{\rm deg\, \Delta_{11}}}\Delta_{11}=(27+\si^3)^{1/6}x_0\pi_A,\\
&&(-1)^{\frac{1}{2}-{\rm deg\, \Delta_{12}}}\Delta_{12}=\Big(\eta x_0^2 +\frac{1}{3}\Big(\eta \si + (27+\si^3)^{1/3}\Big)x_1\Big)\pi_A, \\
&&(-1)^{\frac{1}{2}-{\rm deg\, \Delta_{13}}}\Delta_{13} = \frac{1}{9}(27+\si^3)^{1/2}\Big(\Phi_{52}(\si)\, x_0^3 -\Phi_{42}(\si)\, x_0x_1\Big)\pi_A,\\
&&(-1)^{\frac{1}{2}-{\rm deg\, \Delta_{14}}}\Delta_{14}=\Big(\frac{1}{3\si^2}
\Big(\eta^2 (18+\si^3)-\si (27+\si^3)^{2/3} \Big) x_1^2
-\eta^2 x_0^4/\si \Big)\pi_A, \\
&& (-1)^{\frac{1}{2}-{\rm deg\, \Delta_{15}}}\Delta_{15}=\frac{(27+\si^3)^{5/6}}{24\si+\si^4}x_0^5\pi_A,\\
&&(-1)^{\frac{1}{2}-{\rm deg\, \Delta_{21}}}\Delta_{21}=
\Big(\frac{1}{3}\Big(-\eta\si+2(27+\si^3)^{1/3}\Big)x_1-\eta x_0^2\Big)\pi_A,\\
&&
(-1)^{\frac{1}{2}-{\rm deg\, \Delta_{22}}}\Delta_{22}=
\Big(-\frac{1}{3\si^2}\Big(\eta^2(18+\si^3)+2\si(27+\si^3)^{2/3}\Big)x_1^2+\eta^2x_0^4/\si\Big)\pi_A,\\
&&(-1)^{\frac{1}{2}-{\rm deg\, \Delta_{31}}}\Delta_{31} = -\frac{1}{2\sqrt{3}}(27+\si^3)^{1/2}\Big(
\Phi_{51}(\si)\, x_0^3 -\Phi_{41}(\si)\, x_0x_1\Big)\pi_A,\\
\een
Since
\ben
{\rm res}_{x=0}\ \frac{x_0x_1x_2}{f_{x_0}f_{x_1}f_{x_2}}\,d^3x =
\frac{3}{2(27+\si^3)}
\een
the identifications for the other two classes should be
\ben
1=1,\quad P=\frac{2}{3}\,(27+\si^3)x_0x_1x_2\,\pi_A^2.
\een
It is easy to check that the Poincar\'e and the residue pairings agree. We also have the analogue of Proposition \ref{p8:prop2}.
\begin{proposition}\label{j10:prop2}
The K\"ahler parameter is related to the modulus $\tau$ via the following mirror map: $t_{-1}=2\pi i\tau/6.$
\end{proposition}
The proof is again along the same lines and it is omitted. In order to recall the main result of Krawitz--Shen we just need to check that the Fourier expansions of the correlators in powers of $q:=e^{2\pi i\tau/6}$ have the correct leading terms. For the first correlator we have
\ben
\langle\Delta_{11},\Delta_{11},\Delta_{14}\rangle_{0,3}=
\frac{1}{18}\Big(\si+2\eta^2(27+\si^3)^{1/3}\Big)\, (-1)^{1/2} \, \pi_A
\een
and the Fourier expansion is the following:
\ben
\langle\Delta_{11},\Delta_{11},\Delta_{14}\rangle_{0,3}=\frac{1}{6}+q^6+q^{18}+q^{24}+O(q^{30}).
\een
The remaining correlators can be computed similarly. The computations are straightforward but quite cumbersome. Here is what we got (with the help of computer software):
\ben
&&\langle\Delta_{11},\Delta_{12},\Delta_{13}\rangle_{0,3}=\frac{1}{6}+q^{12}+O(q^{31}) \\
&&\langle\Delta_{21},\Delta_{21},\Delta_{21}\rangle_{0,3}=\frac{1}{3}+ 2q^6+2q^{18}+2q^{24}+O(q^{30})\\
&&\langle\Delta_{11},\Delta_{21},\Delta_{31}\rangle_{0,3}=q+ 2q^7+2q^{13}+2q^{19}+q^{25}+O(q^{31}).
\een
The main result of Krawitz-Shen in the case of the $J_{10}$ singularity can be formulated in this way:
\begin{theorem}\label{lcsl:j10}
 Let $q=e^{2\pi i\tau/6}$. Under the above identification of the
quantum cohomology and the Milnor ring, $\D_{LCS,q}$ is equal to
the descendant potential function of the elliptic orbifold $\P^1$ with
weights $(6,3,2)$.
 \end{theorem}

Theorems \ref{quasi_mod} and \ref{lcsl:j10} imply

\begin{corollary}
The Gromov--Witten total descendant potential function of the elliptic orbifold $\P^1$ with
weights $(6,3,2)$ is quasi-modular for $q=e^{2\pi i \tau/6}$ and a finite index subgroup
of ${\rm SL}_2(\Z)$.
\end{corollary}

\appendix
\section{Primitive forms for simple elliptic singularities}\label{app:A}

\subsection{Oscillatory integrals}
If $V$ is a vector space, then we denote by $V((z))$ (resp. $V[[z]]$) the space of formal Laurent (resp. power) series in $z$ with coefficients in $V$. The space $\H_F$ of oscillatory integrals is defined formally as the third cohomology of the twisted de Rham complex:
\ben
q_*\Omega^\bullet_{X/S}((z)),\quad d=zd_{X/S} + dF\wedge.
\een
Given a differential form $\omega\in q_*\Omega_{X/S}^3$, we denote by $\int e^{F/z}\omega$ its projection on $\H_F$. The sheaf $\H_F$ is equipped with a Gauss-Manin connection:
\ben
\nabla^{\rm G.M}_{\d/\d t^a} \int e^{F/z}g(t,x,z)d^3x = \int e^{F/z}\Big(z^{-1}\frac{\d F}{\d t^a}g +\frac{\d g}{\d t^a}\Big) d^3x
\een
and
\ben
\nabla^{\rm G.M}_{\d/\d z} \int e^{F/z}g(t,x,z)d^3x = \int e^{F/z}\Big(-z^{-2}Fg +\frac{\d g}{\d z}\Big) d^3x,
\een
where  $d^3x=dx_0dx_1dx_2$ is the standard volume form.

We say that an element $\omega\in \H_F$ is homogeneous of degree $r$ if
\ben
(z\nabla_{\d/\d z} + \nabla_E)\omega = r\omega,
\een
where $E$ is the Euler vector field.
Let us denote by $\H_F^{(0)}$ the subspace of $\H_F$ consisting of power series in $z$. According to K. Saito there exists a sequence of bilinear pairings:
\ben
K_F^{(k)}:\H_F^{(0)}\times \H_F^{(0)}\to \O_\S,\quad k\geq 0,
\een
satisfying the following properties:
\begin{itemize}
\item[(K1)] The pairings are symmetric for $k$ even and skew-symmetric for $k$ odd.
\item[(K2)] The pairings are compatible with the Gauss-Manin connection:
\ben
\xi\, K_F^{(k)}(\omega_1,\omega_2)= K_F^{(k)}(\nabla_\xi\omega_1,\omega_2)+K_F^{(k)}(\omega_1,\nabla_\xi\omega_2)
\een
for all $\xi\in \T_S$ and all $\omega_1,\omega_2\in \H_F^{(0)}.$.
\item[(K3)]
We have $K_F^{(k)}(z\omega_1,\omega_2)=K_F^{(k-1)}(\omega_1,\omega_2).$
\item[(K4)]
If $\omega_1$ and $\omega_2$ are homogeneous of degrees $r_1$ and $r_2$ then $K_F^{(k)}(\omega_1,\omega_2)$ is homogeneous of degree $r_1+r_2-k-3$ (the number 3 here corresponds to the fact that $\omega_1$ and $\omega_2$ are 3-forms).
\item[(K5)]
If $\omega_i\in \H_F$ are represented by differential forms $g_i(t,x)d^3x$ independent of $z$ then
$K_F^{(0)}(\omega_1,\omega_2)$ coincides with the residue pairing
\ben
\Big(\frac{1}{2\pi i}\Big)^{3}\int_{\Gamma_\epsilon}
\frac{g_1(t,x)\ g_2(t,x)}
{F_{x_0}F_{x_1} F_{x_{2}} }\ d^3x\ ,
\een
where the integration cycle $\Gamma_\epsilon$ is supported on
$|\frac{\d F}{\d x_0}|= |\frac{\d F}{\d x_{1}}| =|\frac{\d F}{\d x_{2}}|=\epsilon$.
\end{itemize}

\subsection{The primitive forms}

Let $g(s,x)d^3x$ be a volume form (i.e., $g(s,x)\neq 0$ for all $(s,x)\in \S\times \C^3$) and $\omega\in \H_F^{(0)}$ be the corresponding oscillatory integral.  The period mapping
\beq\label{period_iso}
\d/\d s^a \mapsto z\nabla_{\d/\d s_a}\
\int e^{F/z}g(s,x)d^3x,\ \ \ 1\leq a\leq N,
\eeq
induces an isomorphism between $\T_\S[[z]]$ and $\H_F^{(0)}.$ The volume form is called {\em primitive} if it is homogeneous and it satisfies the following properties:
\begin{enumerate}
\item[(P1)]
For all vector fields $\d/\d s^i$, $\d/\d s^j$ and all $k\geq 1$ we have
\ben
 K_F^{(k)}\Big(z\nabla_{\d/\d s^i} \omega, z\nabla_{\d/\d s^j}\omega\Big)=0
\een
 \item[(P2)]
For all vector fields $\d/\d s^i$, $\d/\d s^j$, $\d/\d s^l$ and all $k\geq 2$ we have
\ben
K_F^{(k)}\Big(z\nabla_{\d/\d s^i}z\nabla_{\d/\d s^j} \omega, z\nabla_{\d/\d s^l}\omega\Big)=0.
\een
\item[(P3)] For all vector fields $\d/\d s^i$, $\d/\d s^j$ and for all $k\geq 2$ we have
\ben
K_F^{(k)}\Big(-z^2\nabla_{\d/\d z}\ z\nabla_{\d/\d s^i} \omega, z\nabla_{\d/\d s^j}\omega\Big)=0.
\een
\end{enumerate}
Since $g(s,x)d^3x$ is a homogeneous volume form the function $g(s,x)$ must be homogeneous of degree $0$, i.e., $g(s,x)$ depends only on the degree-0 variable $s_{-1}=\sigma.$ Condition (P1) holds for any degree-0 function $g$, due to the homogeneity (K4) and the skew-symmetry (K1) of $K_F$.

We are going to prove that all primitive forms have the form $d^3x/\pi_a(\si)$, where $a$ is a monodromy invariant cycle and $\pi_a$ is the Gelfand-Leray period (\ref{lerey_period}).

All identities involving holomorphic forms should be understood in the space $\H_F$ of oscillatory integrals, i.e., we work modulo $(zd_{X/S}+dF\wedge)$-exact forms.
Note that we have the following identity:
\beq\label{relation}
\frac{\d F}{\d s^i}\,\frac{\d F}{\d s^j} d^3x = C_{ij}^k\frac{\d F}{\d s^k} d^3x+ z B_{ij}^k\frac{\d F}{\d s^k}d^3x+z^2A_{ij}^k\frac{\d F}{\d s^k} d^3x,
\eeq
where we adopted Einstein's convention for summation over repeating lower and upper indices. The coefficients are homogeneous functions on $S$ of degree, respectively:
\ben
{\rm deg}\, C_{ij}^k &=& {\rm deg}\, s_k - {\rm deg}\, s_i-{\rm deg}\, s_j+1,  \\
{\rm deg}\, B_{ij}^k &=& {\rm deg}\, s_k -{\rm deg}\,s_i -{\rm deg}\,  s_j ,\\
{\rm deg}\, A_{ij}^k &=& {\rm deg}\, s_k-{\rm deg}\, s_i- {\rm deg}\,   s_j -1.\\
\een
In particular, $A_{ij}^k\neq 0$ only for $i=j=-1$ and $k=0$. Let $\omega=g(\sigma)d^3x\in \H_F$ be a primitive form. A straightforward differentiation gives that $z\nabla_{\d/\d s^i}z\nabla_{\d/\d s^j} \omega$ is a sum of three terms:
\ben
C_{ij}^k(z\nabla_{\d/\d s^k}\omega) ,
\een
\ben
z\,\Big( -C_{ij}^k\frac{1}{g}\frac{\d g}{\d s^k} z\nabla_{\d/\d s^0}\omega + \frac{1}{g}\frac{\d g}{\d \si}
\Big(
\delta_{i,-1}\, z\nabla_{\d/\d s^j}\omega + \delta_{j,-1}\, z\nabla_{\d/\d s^i}\omega \Big)+ B_{ij}^k\, z\, \nabla_{\d/\d s^k}\omega
\Big),
\een
and
\ben
z^2\delta_{i,-1}\delta_{j,-1}\Big(
A_{ij}^0(\si)-\frac{2}{g^2}\Big(\frac{\d g}{\d \si}\Big)^2 - B_{ij}^{-1}\frac{1}{g}\frac{\d g}{\d \si}+\frac{1}{g}\frac{\d^2 g}{\d \si^2}\Big)\, z\nabla_{\d/\d s^0}\omega
\een
In order for $\omega$ to be primitive we have to arrange that
$$
K_F^{(2)}(z\nabla_{\d/\d s^i}z\nabla_{\d /\d s^j}\omega,z\nabla_{\d/\d s^l}\omega)=0\quad \mbox{for all}\  -1\leq i,j,l\leq 6.
$$
We already know that $K_F^{(k)}(z\nabla_{\d/\d s^i}\omega,z\nabla_{\d/\d s^l})=0$ for all $i$ and $l$, and all $k\geq 1$. Therefore, using property (K3) of the higher residue pairing, we get that it is enough to prove that the last of the above 3 terms is 0. In other words, $g$ must be a solution to a second order differential equation. Put $u=1/g$; then the differential equation becomes:
\ben
\frac{\d^2 u}{\d \si^2}= B_{-1,-1}^{-1}(\si)\, \frac{\d u}{\d\si} + A_{-1,-1}^0 (\sigma)\,u.
\een
Comparing with equation (\ref{relation}) we see that the solutions of this differential equation can be constructed via the oscillatory integrals $\int e^{f/z}d^3x$ which are obtained from $\int e^{F/z}d^3x$ by specializing the parameters $s_0=s_1=\dots=s_6=0$. Note that this substitution is necessary in order for $C_{-1,-1}^k$ to become $0$ so that the relation (\ref{relation}) matches the above differential equation.

Alternatively, solutions of the differential equation can be
constructed via the Laplace transform of the oscillatory
integrals. Namely, the Gelfand-Leray periods $\pi_\ga$ where $\ga$ is
any flat middle homology cycles. Note however, that the Gelfand-Leray
periods vanish whenever $\ga$ is an
eigenvector of the classical monodromy with eigenvalue different from
1. Therefore, we may assume that $\ga$ is an invariant cycle with
respect to the classical monodromy, i.e., it is a tube cycle.
\begin{comment}
\section{The monodromy group of $X_9$}\label{app:B}
Note that according to Proposition \ref{x9:prop1} there exists a symplectic basis  $\{A, B\}$ near
$\si=\infty$ such that the local monodromy around $\si=\infty$ is given by
\ben
T_3^{\rm local}=\begin{bmatrix}-1&-2\\0&-1\end{bmatrix} .
\een
Similarly, by matching the solutions of the differential equation \eqref{PF:X9} near the other two
singular points $\si_1:=-2$ and $\si_2:=2$ with the $j$-nvariant (see the proof of Proposition
\ref{p8:prop1}), one can prove that the {\em local} monodromies around $\si=2$ and $\si=-2$ are
\ben
T_1^{\rm local}=T_2^{\rm local}=\begin{bmatrix}1&2\\ 0&1\end{bmatrix}.
\een
Let us fix 3 loops around each
singular point and pick a symplectic basis $\{A,B\}$ such that the
monodromies satisfy $T_2=T_2^{\rm local}$, $T_1=g_1T_1^{\rm  local}g_1^{-1}$, and
$T_3=g_3T_3^{\rm local}g_3^{-1}$ for some matrices $g_1,g_3\in {\rm SL}_2(\Z)$, and
\beq\label{mon:eqn}
T_3T_2T_1=1.
\eeq
Note that $T_1$ and $T_3$ depend on 4 parameters. We find them by solving \eqref{mon:eqn}:
\ben
T_1=
\begin{bmatrix}
1-2p & 2p^2 \\
-2    & 1+2p
\end{bmatrix}
\een
for some integer $p$. Conjugating $T_1$ by $T_2$ decreases/increases $p$ by $2$. Therefore we may
assume that $p=1$ or $p=0$. Finally, changing the symplectic basis $\{A,B\}$ if necessary,
we can arrange that $p=0.$
\end{comment}

\end{document}